\documentclass[10pt,a4paper,reqno]{amsart}

\usepackage[a4paper,margin=31mm]{geometry}
\usepackage[T1]{fontenc}
\usepackage[utf8]{inputenc}
\usepackage{lmodern}
\usepackage{microtype}
\usepackage{amsmath,amssymb,amsthm,mathtools}
\usepackage{enumitem}
\usepackage{booktabs,tabularx,array}
\usepackage{xcolor}
\usepackage{aliascnt}
\usepackage{hyphenat}
\usepackage{verbatim}
\usepackage[hidelinks]{hyperref}
\usepackage[nameinlink,noabbrev]{cleveref}

\allowdisplaybreaks
\setlist[itemize]{leftmargin=2em,itemsep=2pt,topsep=4pt}
\setlist[enumerate]{leftmargin=2.3em,itemsep=2pt,topsep=4pt}
\newtheorem{theorem}{Theorem}[section]

\newaliascnt{proposition}{theorem}
\newtheorem{proposition}[proposition]{Proposition}
\aliascntresetthe{proposition}

\newaliascnt{lemma}{theorem}
\newtheorem{lemma}[lemma]{Lemma}
\aliascntresetthe{lemma}

\newaliascnt{corollary}{theorem}
\newtheorem{corollary}[corollary]{Corollary}
\aliascntresetthe{corollary}

\newaliascnt{claim}{theorem}

\aliascntresetthe{claim}

\theoremstyle{definition}
\newaliascnt{definition}{theorem}
\newtheorem{definition}[definition]{Definition}
\aliascntresetthe{definition}

\theoremstyle{remark}
\newaliascnt{remark}{theorem}
\newtheorem{remark}[remark]{Remark}
\aliascntresetthe{remark}

\crefname{theorem}{theorem}{theorems}
\Crefname{theorem}{Theorem}{Theorems}
\crefname{proposition}{proposition}{propositions}
\Crefname{proposition}{Proposition}{Propositions}
\crefname{lemma}{lemma}{lemmas}
\Crefname{lemma}{Lemma}{Lemmas}
\crefname{corollary}{corollary}{corollaries}
\Crefname{corollary}{Corollary}{Corollaries}
\crefname{remark}{remark}{remarks}
\Crefname{remark}{Remark}{Remarks}

\newcommand{\C}{\mathbb C}
\newcommand{\R}{\mathbb R}

\newcommand{\Z}{\mathbb Z}
\newcommand{\GL}{\operatorname{GL}}
\newcommand{\SL}{\operatorname{SL}}
\newcommand{\PSL}{\operatorname{PSL}}
\newcommand{\Mat}{\operatorname{Mat}}
\newcommand{\Ad}{\operatorname{Ad}}
\newcommand{\ad}{\operatorname{ad}}
\newcommand{\tr}{\operatorname{tr}}
\newcommand{\id}{\operatorname{id}}

\newcommand{\Hom}{\operatorname{Hom}}

\newcommand{\K}{\mathbb K}

\newcommand{\cG}{\mathcal G}
\newcommand{\cT}{\mathcal T}

\newcommand{\cE}{\mathcal E}
\newcommand{\abs}[1]{\lvert#1\rvert}

\title{$3$-manifold groups are linear}

\author{Xiaolei Wu}
\address{Shanghai Center for Mathematical Sciences, Jiangwan Campus, Fudan University, No. 2005 Songhu Road, Shanghai, 200438, P.R. China}
\email{xiaoleiwu@fudan.edu.cn}

\author{Shengkui Ye}
\address{NYU Shanghai, No. 567 Yangsi West 
  Rd, Pudong New Area, Shanghai, 200124, P.R. China \\
NYU-ECNU Institute of Mathematical Sciences at NYU Shanghai, 3663 Zhongshan Road North, Shanghai, 200062, China}
\email{sy55@nyu.edu}

\subjclass[2020]{57M05, 20H20, 20E08}
\keywords{$3$-manifold groups, linear groups, graphs of groups}
\date{October 2026}

\begin{document}
\begin{abstract}
We prove that the fundamental group of every compact graph $3$-manifold is
linear over the complex numbers.  As a corollary, all  fundamental groups of compact  $3$-manifolds are linear over the complex numbers.
\end{abstract}

\maketitle


\section*{Introduction}
\label{sec:introduction}

A group $G$ is \emph{linear over a field $k$} if it admits a faithful
finite-dimensional representation, that is, an injective homomorphism
\[
  G \longrightarrow \GL_n(k)
\]
for some positive integer $n$.  Linearity  imposes strong restrictions on the group
structure. By a theorem of Mal'cev, every finitely generated linear group
is residually finite \cite{Mal40}. In characteristic zero, Selberg's lemma
gives a torsion-free subgroup of finite index \cite{Sel60}, while the Tits
alternative asserts that every finitely generated subgroup is either
virtually solvable or contains a nonabelian free group \cite{Tit72}.

Many important classes of groups are known to be linear. Classical
examples include finitely generated free groups, fundamental groups of
compact surfaces, virtually polycyclic groups \cite{Swn67}, and finitely
generated Coxeter groups. Bigelow and
Krammer independently proved that braid groups are linear
\cite{Big01,Kra02}. Davis and Januszkiewicz established linearity of
finitely generated right-angled Artin groups by showing that they are commensurable with
right-angled Coxeter groups \cite{DavisJanuszkiewicz}. The theory of special cube
complexes developed by Haglund and Wise provides a further source of
examples: finitely generated virtually special groups virtually embed
in right-angled Artin groups and are therefore linear \cite{HW08}.
Linearity over the complex numbers is preserved by taking subgroups, finite-index
overgroups, and finite direct products. It is also preserved by finite
graph products, as proved by Berlai and de la Nuez \cite{BN22}; this
includes finite free products.
For amalgamated products and HNN extensions of linear groups, one usually needs to require strong hypotheses to make them linear; see
\cite{Wehrfritz73,Shalen79,DoubaTs} for some progress in this direction.

Whether $3$-manifold groups are linear has been a longstanding question
in low-dimensional topology. Here by $3$-manifold group we mean the fundamental group of a compact $3$-manifold. In Kirby's problem list, Thurston asked
whether every  $3$-manifold group admits a faithful
representation into $\GL_4(\R)$
\cite[Problem~3.33(A)]{Kirby}. Button showed that this proposed dimension
bound is impossible: some closed orientable graph-manifold groups do not
embed in $\GL_4(k)$ over any field $k$ \cite{Button}. The question of linearity without restrictions on the representation
dimension remained open and appears in many places; see, for example, \cite[Question 7.3.1 (2)]{AFW}, \cite[Question 2]{AgolICM}, 
\cite[Question~7.6]{Boileau} and the new Kirby problem list \cite[Problem~3.37(a)]{K3}.

Geometrization and the theory of virtually special groups settled most
cases of the linearity question. Passing, if necessary, to the orientable double cover and
using prime decomposition,   the question is reduced to compact orientable irreducible $3$-manifolds. Hyperbolic
$3$-manifold groups are linear through their holonomy representations in
$\PSL_2(\C)$, and Seifert and Sol groups are also known to
be linear \cite{AFW}.  Przytycki and Wise proved that
mixed $3$-manifold groups are virtually special \cite{PWMixed}. Here a
mixed manifold is a compact orientable irreducible $3$-manifold that is
neither hyperbolic nor a graph manifold. These results reduce the linearity problem for $3$-manifold groups to graph manifold case.

For graph manifolds, whose JSJ pieces are Seifert fibered, Liu proved
that an aspherical compact graph manifold admits a nonpositively curved
metric if and only if it is virtually special \cite{Liu}. Przytycki and Wise proved that
graph manifolds with nonempty boundary are virtually special
\cite{PWBoundary}. Consequently, these groups are linear. The remaining
case consists of closed nongeometric graph manifolds, i.e., they do not admit
nonpositively curved metrics. Typical examples are mapping tori of closed surfaces of genus at least two, with monodromy given by a Dehn twist along an essential simple closed curve.

Residual properties provide a related perspective on the problem. A
group is \emph{residually finite} if every nontrivial element has a
nontrivial image in some finite quotient. For a prime $p$, it is
\emph{residually $p$} if the quotient can always be chosen to be a finite
$p$-group, and \emph{virtually residually $p$} if it has a finite-index
subgroup with this property. Hempel's work, together with geometrization,
implies that fundamental groups of compact $3$-manifolds are residually
finite \cite{Hem87,AFW}. Aschenbrenner and Friedl proved that every
such group is virtually residually $p$ for all but finitely many primes,
and for every prime in the graph-manifold and fibered cases
\cite{AF13}. Combined with the subsequent virtual specialness and
virtual fibering results, their work gives virtual residual $p$ for
every prime $p$ for all compact $3$-manifold groups; see
\cite[Introduction]{Wil17}. These results support but do not prove  that $3$-manifold groups are linear. 

The purpose of this paper is to settle the remaining graph-manifold
case.

\begin{theorem}[Main theorem]\label{thm:main}
Let $M$ be a closed graph $3$-manifold. Then there are an integer
$N\geq 1$ and an injective homomorphism
\[
  \pi_1(M)\hookrightarrow \GL_N(\C).
\]
In fact, when $M$ is  nongeometric, it has a finite cover  whose fundamental group admits a faithful  representation in dimension $73$.
\end{theorem}

Together with the known cases recalled above, this gives the following
consequence.

\begin{corollary}\label{cor:3-mfd-linear}
The fundamental group of every compact $3$-manifold is linear over
$\C$.
\end{corollary}

\begin{remark}
   In \cite[Problem~3.37(b)]{K3}, it was further asked whether all $3$-manifold groups could have faithful linear representations with a uniform bound on the dimension. Our construction does not provide such a bound. One reason is that in \Cref{thm:main}, we do not have  good control over the index of the cover.
\end{remark}

Our proof constructs faithful representations directly from the JSJ
decomposition. For the remaining closed orientable irreducible
nongeometric graph manifolds, we first pass to a finite cover whose
Seifert blocks are products $Q_v\times S^1$, where $Q_v$ is a compact
orientable surface with nonempty boundary and genus at least two
\cite{WY,KapovichLeeb}.
The resulting graph of groups has vertex groups
$\pi_1(Q_v)\times\Z$ and edge groups $\Z^2$. Although
the vertex groups are linear, their representations must be chosen to
agree on the peripheral subgroups under the prescribed torus gluings.
One must then ensure that the resulting representation is faithful on
elements traversing several JSJ blocks.

We address these two requirements by establishing a new combination theorem for graphs of groups that may be of independent interest. The argument uses Laurent-polynomial deformations
of the edge intertwiners. For a cyclically reduced Bass--Serre word of
positive edge length, the transition conditions ensure that its matrix representative has
a nonzero coefficient at a nonconstant Laurent monomial.
A rank-one construction reduces these conditions to vectors and
covectors detecting the incident edge subgroups. For the graph-manifold
application, we construct vertex representations with
peripheral detectors that are compatible with arbitrary torus gluing
matrices. This yields a faithful representation of the product-block
cover in $\GL_{73}(\R)$.

\subsection*{AI tool disclosure} Our proof is based on an $7$-page draft generated by ChatGPT 5.6 Sol in July 2026. All mathematical content was carefully reviewed, verified, and substantially revised by the authors, who take full responsibility for the final text.

\subsection*{Acknowledgments}
Wu is currently a member of LMNS at Fudan University. After we finished our draft and started circulating it among a few colleagues, we learned that Daniel Wise has independently established a proof of Theorem \ref{thm:main} \cite{Wise26}.

\section{A combination theorem}
\label{sec:combination}

We begin by fixing the graph-of-groups conventions and recalling the
normal-form facts used in the faithfulness argument.

\subsection*{Graphs of groups and their fundamental groups}

Let $\Gamma$ be a finite connected graph, possibly with loops and
multiple edges. Write $V(\Gamma)$ and $E^{\mathrm{or}}(\Gamma)$ for its
vertices and oriented edges. For an oriented edge $e$, let $o(e)$ and
$t(e)$ be its initial and terminal vertices, and let $\bar e$ be the
same edge with the opposite orientation.

A \emph{graph of groups} $\mathcal G$ over $\Gamma$ consists of a group
$G_v$ for each vertex $v$, a group $G_e=G_{\bar e}$ for each geometric
edge $|e|$, and an injective homomorphism
\[
\iota_e:G_e\longrightarrow G_{o(e)}
\]
for each oriented edge $e$. Thus the two maps $\iota_e$ and
$\iota_{\bar e}$ attach the edge group to its two endpoint groups.
We write
\[
H_e:=\iota_e(G_e)\leq G_{o(e)}.
\]

An \emph{incidence} is an oriented edge marked at its initial vertex.
Note that incidences are in one-to-one correspondence with oriented edges, but it is sometimes more convenient to talk about incidences. 
Each geometric edge therefore has two incidences, including when it
is a loop. Let $\mathcal E$ be the set of geometric edges, and write
$|e|\in\mathcal E$ for the edge underlying $e$; in particular,
$|e|=|\bar e|$.

Fix a base vertex $v_*$ and a maximal subtree $T\subseteq\Gamma$.
The \emph{fundamental group} $\pi_1(\mathcal G,v_*,T)$ is obtained
from the free product
\[
\left(\mathop{*}_{v\in V(\Gamma)}G_v\right)
*
F\bigl(t_e:e\in E^{\mathrm{or}}(\Gamma)\bigr),
\]
where $F$ denotes the free group on the indicated symbols, by
imposing the relations
\[
\begin{aligned}
t_{\bar e}&=t_e^{-1}
&&\bigl(e\in E^{\mathrm{or}}(\Gamma)\bigr),\\
t_e&=1
&&\bigl(|e|\in E(T)\bigr),\\
t_e\iota_e(h)t_e^{-1}&=\iota_{\bar e}(h)
&&\bigl(e\in E^{\mathrm{or}}(\Gamma),\ h\in G_e\bigr).
\end{aligned}
\]
Its isomorphism type is independent of the chosen base vertex and
maximal subtree. The normal-form theorem implies that each natural
map $G_v\to\pi_1(\mathcal G,v_*,T)$ is injective, so we regard the
vertex groups as subgroups of the fundamental group.

\subsection*{Reduced paths and cyclic normal form}

The edge relations allow a traversal of $e$ followed by $\bar e$
to be cancelled when the intervening vertex-group element belongs
to $H_{\bar e}$. This gives the following notion of a reduced
junction.

\begin{definition}[Reduced junctions]
\label{def:reduced-junction}
Let $e,f$ be oriented edges with $t(e)=o(f)=v$, and let $g\in G_v$.
The triple $(e,g,f)$ is \emph{reduced} if
\[
f\neq\bar e,
\qquad\text{or}\qquad
f=\bar e\ \text{ and }\ g\notin H_{\bar e}.
\]
\end{definition}

For the cyclic normal form, we use closed Bass--Serre paths written
as
\[
c=e_1g_1e_2g_2\cdots e_mg_m,
\qquad m\geq1,
\]
where there are vertices $v_0,\ldots,v_m$ satisfying
\[
o(e_j)=v_{j-1},\qquad
t(e_j)=v_j,\qquad
g_j\in G_{v_j},
\qquad v_m=v_0.
\]
The expression records the order of traversal: first $e_1$, then
$g_1$ at $v_1$, then $e_2$, and so on. We write group products in
composition order, with the rightmost factor acting first.
Accordingly, the element associated with $c$ is
\[
\gamma_c=g_mt_{e_m}\cdots g_2t_{e_2}g_1t_{e_1}.
\]
Note that here $v_0$ might not equal to $ v_*$, but we can always join $v_*$ to $v_0$ along the unique path in
$T$ and return along its reverse. Since the edge letters in $T$ are
trivial, this gives the same formula for $\gamma_c$. 

The path $c$ is \emph{reduced} if each junction
$(e_j,g_j,e_{j+1})$ is reduced for $1\leq j<m$. It is
\emph{cyclically reduced} if the closing junction
$(e_m,g_m,e_1)$ is reduced as well. Equivalently, setting
$e_{m+1}=e_1$, cyclic reducedness means
\[
e_{j+1}=\bar e_j
\quad\Longrightarrow\quad
g_j\notin H_{\bar e_j}
\qquad(1\leq j\leq m).
\]

The normal-form theorem says that every reduced closed path of
positive edge length represents a nontrivial element; see
\cite[Chapter~I, \S5.2, Theorem~11 and Corollary~3]{Serre}.
A standard consequence is the following: 
\begin{theorem}[Cyclic normal form theorem]\label{thm-cyc-nor-form}
    Every
nontrivial element of $\pi_1(\mathcal G,v_*,T)$ is conjugate either
to a nontrivial element of a vertex group or to $\gamma_c$ for a
cyclically reduced closed path $c$ of positive edge length.
\end{theorem}
\begin{proof}[Sketch of proof]
To obtain such a representative, one first reduces it to a closed word.
If its closing junction is not reduced, a cyclic permutation,
which corresponds to conjugation, makes that junction an internal
one. Cancelling the resulting inverse edges and absorbing the
vertex-group factors gives a shorter closed word. Repeating this
procedure terminates at a vertex-group element or a cyclically
reduced closed path.
\end{proof}

\subsection*{Leading-block combination theorem}
We now formulate the representation data that ensure these reduced
paths have nontrivial matrix images.

\begin{definition}[Leading-block system of rank $d$]\label{defn-leading-block}
Let $k$ be a field.  A \emph{leading-block system of rank $d$} for
$\cG$ consists of the following data.

\begin{enumerate}[label=\rm(\roman*)]

\item
For every vertex $v$, there is an $n$-dimensional $k$-vector space $V_v$ and a
faithful representation
\[
\sigma_v\colon G_v\longrightarrow \GL(V_v).
\]

\item
For every geometric edge $|e|\in\cE$, there is an indeterminate $x_{|e|}$, such that the family
$\{x_{|e|}:|e|\in\cE\}$ is algebraically independent over $k$.  Set
\[
K=k(x_{|e|}:|e|\in\cE).
\]

\item
For every oriented edge $e$, there is a $K$-linear isomorphism
\[
B_e\colon V_{o(e)}\otimes_k K
\longrightarrow
V_{t(e)}\otimes_k K
\]
such that
\[
B_{\bar e}=B_e^{-1}
\]
and
\[
B_e\sigma_{o(e)}\bigl(\iota_e(h)\bigr)
=
\sigma_{t(e)}\bigl(\iota_{\bar e}(h)\bigr)B_e
\qquad (h\in G_e).
\]

\item
There is a fixed $d$-dimensional $k$-vector space $W$ such that, for
every oriented edge $e$,
\[
B_e\in
\Hom_k\bigl(V_{o(e)},V_{t(e)}\bigr)
   \otimes_k k[x_{|e|}^{\pm1}],
\qquad
\deg_{x_{|e|}}B_e\leq 1,
\]
and there are $k$-linear maps
\[
p_e\colon V_{o(e)}\longrightarrow W,
\qquad
i_e\colon W\longrightarrow V_{t(e)}
\]
satisfying
\[
[x_{|e|}]B_e=i_e\circ p_e.
\]
Here $[x_{|e|}]B_e$ denotes the coefficient of $x_{|e|}^{\,1}$ in
the Laurent-polynomial map $B_e$.

\item
For every reduced triple $(e,g,f)$, with $t(e)=o(f)=v$, 
the transition map
\[
p_f\sigma_v(g)i_e\colon W\longrightarrow W
\]
is invertible.

\end{enumerate}
\end{definition}

\begin{theorem}[Leading-block combination theorem]
\label{thm:leading-block}
If $\mathcal G$ admits a leading-block system of rank $d$ for some $d\geq 1$, then $\pi_1(\mathcal G,v_*,T)$ embeds in $\GL_n(K)$.

\end{theorem}

\begin{proof}
We extend each representation $\sigma_v$ from $V_v$ to
$V_v\otimes_k K$ by extension of scalars, without changing the
notation.

For every vertex $v$, let $\tau_v$ be the unique path in the maximal
subtree $T$ from $v_*$ to $v$.  If $\tau_v$ traverses the oriented
edges
\[
a_1,\ldots,a_s
\]
in this order, define
\[
U_v=B_{a_s}\cdots B_{a_1}
  \colon V_{v_*}\otimes_k K
  \longrightarrow V_v\otimes_k K,
\]
and put
\[
U_{v_*}=I.
\]
Every $U_v$ is a $K$-linear isomorphism.

For $g\in G_v$, define
\[
\rho(g)=U_v^{-1}\sigma_v(g)U_v,
\]
and, for every oriented edge $e$, define
\[
\rho(t_e)=U_{t(e)}^{-1}B_eU_{o(e)}.
\]
Thus all the matrices $\rho(g)$ and $\rho(t_e)$ belong to
\[
\GL(V_{v_*}\otimes_kK)\cong\GL_n(K).
\]

We first verify the Bass--Serre relations.  Since
$B_{\bar e}=B_e^{-1}$, we have
\[
\begin{aligned}
\rho(t_{\bar e})
&=
U_{o(e)}^{-1}B_{\bar e}U_{t(e)}\\
&=
U_{o(e)}^{-1}B_e^{-1}U_{t(e)}
=
\rho(t_e)^{-1}.
\end{aligned}
\]

Suppose that the geometric edge underlying $e$ belongs to $T$.  Then
\[
U_{t(e)}=B_eU_{o(e)}.
\]
Indeed, if the path from $v_*$ to $t(e)$ is obtained by appending $e$
to the path from $v_*$ to $o(e)$, this follows immediately from the
definition of $U_v$.  In the opposite case, the path from $v_*$ to
$o(e)$ is obtained by appending $\bar e$ to the path from $v_*$ to
$t(e)$, and hence
\[
U_{o(e)}
=
B_{\bar e}U_{t(e)}
=
B_e^{-1}U_{t(e)}.
\]
This gives the same identity.  Consequently,
\[
\rho(t_e)=I
\qquad (e\subset T).
\]

Finally, let $h\in G_e$.  Using the intertwining relation in
Definition~\ref{defn-leading-block}, we obtain
\[
\begin{aligned}
\rho(t_e)\rho\bigl(\iota_e(h)\bigr)\rho(t_e)^{-1}
&=
U_{t(e)}^{-1}
B_e\sigma_{o(e)}\bigl(\iota_e(h)\bigr)
B_e^{-1}U_{t(e)}\\
&=
U_{t(e)}^{-1}
\sigma_{t(e)}\bigl(\iota_{\bar e}(h)\bigr)
U_{t(e)}\\
&=
\rho\bigl(\iota_{\bar e}(h)\bigr).
\end{aligned}
\]
Thus the defining Bass--Serre relations are satisfied.  Therefore
the preceding assignments determine a homomorphism
\[
\rho\colon
\pi_1(\mathcal G,v_*,T)
\longrightarrow
\GL(V_{v_*}\otimes_k K)
\cong\GL_n(K).
\]

We now prove that $\rho$ is injective.  By the cyclic normal-form
theorem for graphs of groups (Theorem \ref{thm-cyc-nor-form}), every nontrivial element is conjugate
either to a nontrivial element of some vertex group or to an element
represented by a cyclically reduced closed Bass--Serre path of
positive edge length.  Since membership in $\ker\rho$ is invariant
under conjugation, it suffices to consider these two cases.

First suppose that the cyclically reduced representative has positive
edge length $m$.  Choose its initial vertex $v_0$ immediately before
an edge and write the corresponding closed path as
\[
c=e_1g_1e_2g_2\cdots e_mg_m.
\]
Put
\[
o(e_j)=v_{j-1},
\qquad
t(e_j)=v_j,
\qquad
g_j\in G_{v_j},
\]
where
\[
v_m=v_0.
\]
We read the indices cyclically and set $e_{m+1}=e_1$.  Cyclic
reducedness means that every cyclic junction
\[
e_j\,g_j\,e_{j+1}
\]
is reduced.  In particular,
\[
e_{j+1}=\bar e_j
\quad\Longrightarrow\quad
g_j\notin H_{\bar e_j}.
\]

Note that the expression
\[
c=e_1g_1e_2g_2\cdots e_mg_m
\]
records the path in chronological order: one first traverses $e_1$,
then applies $g_1$ at $t(e_1)$, then traverses $e_2$, and so on.
Since our matrices act on column vectors on the left, compositions
are written from right to left.  Thus the group element corresponding
to $c$ is
\[
\gamma_c
=
g_mt_{e_m}\cdots g_2t_{e_2}g_1t_{e_1}.
\]
This does not mean that the path is traversed backward: the rightmost
factor $t_{e_1}$ acts first, followed by $g_1$, then $t_{e_2}$, and
so forth.

Substituting the definitions of $\rho(g_j)$ and $\rho(t_{e_j})$,
all intermediate factors $U_{v_j}^{\pm1}$ cancel.  Therefore
\[
\rho(\gamma_c)
=
U_{v_0}^{-1}A_cU_{v_0},
\]
where
\[
A_c
=
\sigma_{v_m}(g_m)B_{e_m}
\sigma_{v_{m-1}}(g_{m-1})B_{e_{m-1}}
\cdots
\sigma_{v_1}(g_1)B_{e_1}.
\]
It is therefore enough to prove that
\[
A_c\neq I.
\]
This is why we compute the leading coefficient of $A_c$, rather than
directly computing a coefficient of the conjugate
$\rho(\gamma_c)$.

Fix a geometric edge $|e|\in\cE$.  For every $j$ such that
$|e_j|=|e|$, write
\[
B_{e_j}=\sum_{r\leq1}x_{|e|}^rB_{j,r}.
\]
Consider one term in the full expansion of the product defining $A_c$.
For every $j$ with $|e_j|=|e|$, suppose that the term
\[
x_{|e|}^{r_j}B_{j,r_j}
\]
is selected from the factor $B_{e_j}$.  Thus $r_j$ is the exponent
selected from the $j$th edge factor, and $r_j\leq1$.  Since the vertex
matrices contain no Laurent variables, the exponent of $x_{|e|}$ in the
resulting term is
\[
\sum_{\{j:|e_j|=|e|\}}r_j.
\]
Define $m_{|e|}:=\#\{j:|e_j|=|e|\}$, then there are $m_{|e|}$ summands in this expression, each at most $1$, so the
exponent is at most $m_{|e|}$.  Equality holds only when $r_j=1$ for every
$j$ satisfying $|e_j|=|e|$.

Consequently,
\[
M_c
=
\prod_{e\in\cE}x_{|e|}^{m_{|e|}}
=
\prod_{j=1}^m x_{|e_j|}
\]
is componentwise maximal among the Laurent monomials occurring in
$A_c$.  Moreover, equality in every coordinate is possible only by taking
the coefficient of $x_{|e_j|}^{\,1}$ from every factor $B_{e_j}$.
Consequently, the coefficient of $M_c$ in $A_c$ is

\[
\begin{aligned}
C_c:=[M_c]A_c
={}&
\sigma_{v_m}(g_m)(i_{e_m}p_{e_m})
\sigma_{v_{m-1}}(g_{m-1})(i_{e_{m-1}}p_{e_{m-1}})
\\
&\qquad\cdots
\sigma_{v_1}(g_1)(i_{e_1}p_{e_1}) ~\\
={}&
\sigma_{v_m}(g_m)i_{e_m} 
\bigl(
p_{e_m}\sigma_{v_{m-1}}(g_{m-1})i_{e_{m-1}}
\bigr)\\
&\quad\cdots
\bigl(
p_{e_2}\sigma_{v_1}(g_1)i_{e_1}
\bigr)
p_{e_1},
\end{aligned}
\]
where $[M_c]A_c$ means that the coefficient of $M_c$ is taken in
every matrix entry.  When $m=1$, the product of the displayed middle
factors is understood to be the identity map on $W$.

For $1\leq j\leq m$, with cyclic indices, define
\[
T_j
=
p_{e_{j+1}}\sigma_{v_j}(g_j)i_{e_j}
\colon W\longrightarrow W.
\]
Since every junction $e_jg_je_{j+1}$ is reduced,
Definition~\ref{defn-leading-block}\textup{(v)} implies that every $T_j$ is invertible.
This includes the last-to-first transition
\[
T_m
=
p_{e_1}\sigma_{v_m}(g_m)i_{e_m},
\]
corresponding to the cyclic junction
\[
e_m\,g_m\,e_1.
\]

Using these transition maps, the leading coefficient becomes
\[
C_c
=
\sigma_{v_m}(g_m)i_{e_m}
T_{m-1}\cdots T_1p_{e_1}.
\]
Since
\[
T_m
=
p_{e_1}\sigma_{v_m}(g_m)i_{e_m}
\]
is invertible, the map $i_{e_m}$ is injective and the map $p_{e_1}$
is surjective.  Hence
\[
\sigma_{v_m}(g_m)i_{e_m}
\]
is injective, while $T_{m-1}\cdots T_1$ is invertible.  It follows
that
\[
\operatorname{rank}C_c=d.
\]
In particular,
\[
[M_c]A_c=C_c\neq0.
\]

Because $m>0$, the monomial $M_c$ is nonconstant.  The identity map
has zero coefficient at every nonconstant Laurent monomial, whereas
$A_c$ has the nonzero coefficient $C_c$ at $M_c$.  Therefore
\[
A_c\neq I,
\]
and consequently
\[
\rho(\gamma_c)\neq I.
\]

It remains to consider a cyclically reduced representative of edge
length zero.  Such an element is a nontrivial element $g\in G_v$ for
some vertex $v$.  Since $\sigma_v$ is faithful,
\[
\sigma_v(g)\neq I,
\]
and therefore
\[
\rho(g)
=
U_v^{-1}\sigma_v(g)U_v
\neq I.
\]
Thus every nontrivial element has nontrivial image under $\rho$.
Hence $\ker\rho=1$, and $\rho$ is faithful.
\end{proof}

\begin{remark}
The rank-one case is particularly transparent.  If
\[
[x_{|e|}]B_e=u_e\ell_e
\]
with $u_e$ a column vector and $\ell_e$ a row vector, then the coefficient of the
maximal monomial is obtained by multiplying the scalar transition coefficients
\[
\ell_f\sigma_v(g)u_e.
\]
The theorem only requires these scalars to be nonzero.
\end{remark}

\section{Rank-one construction} \label{sec:rank-one}

In this section, we specialize the leading-block combination theorem
(\Cref{thm:leading-block}) to the rank-one case and formulate a criterion
in terms of vectors and covectors in the vertex representations.
The resulting theorem (\Cref{prop:combination}) will be applied to
graph $3$-manifold groups in \S~\ref{sec:product-blocks}.

Let $\Gamma$ be a finite connected graph,  $\mathcal G$ be a graph of groups over $\Gamma$ and $k$ be a feild. Suppose that every vertex group $G_v$ has a faithful representation
\[
\sigma_v:G_v\longrightarrow\GL(V_v),
\qquad \dim_kV_v=n.
\]
For every incidence $a$ marked at $v$, suppose there are a vector
\[
Z_a^0\in V_v
\]
and a covector
\[
\lambda_a^0\in V_v^*=\Hom_k(V_v,k)
\]
such that
\[
\lambda_a^0(Z_a^0)=1,
\]
the edge group $H_a$ fixes $Z_a^0$,
\[
\sigma_v(h)Z_a^0=Z_a^0,
\qquad h\in H_a,
\]
and preserves $\lambda_a^0$,
\[
\lambda_a^0\sigma_v(h)=\lambda_a^0,
\qquad h\in H_a.
\]
Finally, assume the pointed detector property
\begin{equation} \label{eq:pi-dec-prop}
\lambda_a^0\sigma_v(g)Z_a^0=1
\quad\Longleftrightarrow\quad
g\in H_a.
\end{equation}

For each incidence $a$ introduce an algebraically independent variable $s_a$ and put
\[
Z_a=s_aZ_a^0,
\qquad
\lambda_a=s_a^{-1}\lambda_a^0.
\]
Then $\lambda_a(Z_a)=1$, and the pointed detector property is unchanged.

For an incidence $a$, let
\[
P_a^0=Z_a^0\lambda_a^0.
\]
Since $\lambda_a^0(Z_a^0)=1$, this is an idempotent.  Define
\[
R_a=s_aP_a^0+I-P_a^0.
\]
Then
\[
R_a^{-1}=s_a^{-1}P_a^0+I-P_a^0.
\]

Suppose that for every preferred orientation $e$ there is an isomorphism
\[
Q_e^0:V_{o(e)}\longrightarrow V_{t(e)}
\]
satisfying the edge-group intertwining relation: for every $h\in G_{\abs e}$,
\[
 Q_e^0\sigma_{o(e)}(\iota_e(h))
 =\sigma_{t(e)}(\iota_{\bar e}(h))Q_e^0,
\]
and
\[
Q_e^0Z_e^0=Z_{\bar e}^0,
\qquad
\lambda_{\bar e}^0Q_e^0=\lambda_e^0.
\]
 Now define
\[
 Q_{\bar{e}}^0 = (Q_e^0)^{-1} , Q_e=R_{\bar e}Q_e^0R_e^{-1},
\]
for each oriented edge $e$.

\begin{lemma}
\label{lem:reverse}
For every oriented edge $e$,
\[
Q_{\bar e}=Q_e^{-1},
\]
and
\[
Q_eZ_e=Z_{\bar e},
\qquad
\lambda_{\bar e}Q_e=\lambda_e.
\]
Moreover,
\[
Q_e\sigma_{o(e)}(\iota_e(h))
=
\sigma_{t(e)}(\iota_{\bar e}(h))Q_e
\qquad(h\in G_e).
\]
\end{lemma}

\begin{proof}
The edge group fixes $Z_e^0$ and preserves $\lambda_e^0$, hence it commutes with
$P_e^0$ and $R_e$.  Therefore the intertwining relation for $Q_e^0$ gives the
corresponding relation for $Q_e$.

Furthermore,
\[
R_e^{-1}Z_e
=
R_e^{-1}(s_eZ_e^0)
=
Z_e^0,
\]
and hence
\[
Q_eZ_e
=
R_{\bar e}Q_e^0Z_e^0
=
R_{\bar e}Z_{\bar e}^0
=
Z_{\bar e}.
\]
The covector identity is similar.

Finally,
\[
Q_{\bar e}
=
R_e(Q_e^0)^{-1}R_{\bar e}^{-1}
=
\left(R_{\bar e}Q_e^0R_e^{-1}\right)^{-1}
=
Q_e^{-1}.
\]
\end{proof}

The point is that $R_e$ and $R_{\bar e}$ need not be inverse to one another.  They
belong to different incidences.  The inverse relation for $Q_e$ follows instead from
the reversal of the order of the two endpoint rescalings.

Enlarge each vertex space by one dimension:
\[
W_v=kb_v\oplus V_v,
\]
where $b_v$ spans the additional one-dimensional summand.  Extend $\sigma_v$ by
\[
D_v(g)=\id \oplus\sigma_v(g).
\]
For every oriented edge $e$, extend $Q_e$ to
\[
\widehat Q_e\colon W_{o(e)}\longrightarrow W_{t(e)}
\]
by
\[
\widehat Q_e(\alpha b_{o(e)}+z)
=
\alpha b_{t(e)}+Q_e(z), 
\qquad \alpha\in k, z\in V_{o(e)}.
\]

Let $b_v^*$ be the coordinate functional defined by
\[
b_v^*(\alpha b_v+z)=\alpha,
\qquad \alpha\in k, z\in V_v,
\]
and extend each $\lambda_a$ at $v$ to $W_v$ by setting
$\lambda_a(b_v)=0$.

For each geometric edge $|e|$, introduce a variable $q_{|e|}$, with the family
$\{q_{|e|}\}$ algebraically independent over the field generated by the matrix
coefficients, the $s_a$, and all previously introduced constants.
We suppress the corresponding extensions of scalars from the notation.
For an incidence $a$ at $v$ with $|a|=|e|$, define
\[
u_a^+=b_v+Z_a,
\qquad
u_a^-=b_v+q_{|e|}Z_a,
\]
and
\[
\ell_a^+
=
\frac{q_{|e|} b_v^*-\lambda_a}{q_{|e|}-1},
\qquad
\ell_a^-
=
\frac{-b_v^*+\lambda_a}{q_{|e|}-1}.
\]
In particular, $q_{|e|}-1\neq0$.

A direct calculation gives
\[
\ell_a^\epsilon(u_a^\delta)=
\begin{cases}
1,&\epsilon=\delta,\\
0,&\epsilon\neq\delta.
\end{cases}
\]
Consequently,
\[
P_a^\pm :=u_a^\pm\ell_a^\pm
\]
are orthogonal rank-one idempotents, i.e.,
\[
(P_a^\pm)^2=P_a^\pm,
\qquad
P_a^+P_a^-=P_a^-P_a^+=0,
\]
and the image of $P_a^\pm$ is of dimension one.

For a variable $x\neq0$, define
\[
C_a(x)
=
xP_a^+
+x^{-1}P_a^-
+I-P_a^+-P_a^-.
\]
Then
\[
C_a(x)^{-1}=C_a(x^{-1}).
\]

Under $D_v$, the subgroup $H_a$ fixes $b_v$ and $Z_a$, and preserves
$b_v^*$ and $\lambda_a$.  Hence $C_a(x)$ commutes with $D_v(H_a)$.

Lemma~\ref{lem:reverse} and the definition of $\widehat Q_e$ give
\[
\widehat Q_{\bar e}=\widehat Q_e^{-1},
\]
and
\[
\widehat Q_eu_e^\pm=u_{\bar e}^\pm,
\qquad
\ell_{\bar e}^\pm\widehat Q_e=\ell_e^\pm.
\]
Therefore
\[
\widehat Q_eC_e(x)\widehat Q_e^{-1}=C_{\bar e}(x).
\]

We  now use these to define the operators $B_e$. Choose one orientation for every geometric edge.  Define
\[
\omega(e)=
\begin{cases}
1,&e\text{ is the chosen orientation},\\
-1,&e\text{ is the reverse orientation},
\end{cases}
\]
and and put
\[
B_e:=\widehat Q_e C_e(x_{|e|}^{\omega(e)})
   :W_{o(e)}\longrightarrow W_{t(e)},
\]
where the vertex spaces are understood to be extended to the
field containing all the parameters $s_a$, $q_{|e|}$, and $x_{|e|}$. Then
\[
B_{\bar e}=B_e^{-1}.
\]
Indeed, since $|\bar e|=|e|$ and $\omega(\bar e)=-\omega(e)$,
\[
B_e^{-1}
=
C_e(x_{|e|}^{-\omega(e)})\widehat Q_e^{-1}
=
\widehat Q_{\bar e}C_{\bar e}(x_{|e|}^{-\omega(e)})
=
B_{\bar e}.
\]
The same identities show that $B_e$ satisfies the required edge-group
intertwining relation.

Let $\epsilon(e)=+$ if $e$ is the chosen orientation; othewise, let
$\epsilon(e)=-$. Since
\[
C_e(x)=xP_e^++x^{-1}P_e^-+I-P_e^+-P_e^-,
\]
we have
\[
[x_{|e|}]\,C_e(x_{|e|}^{\omega(e)})
=P_e^{\epsilon(e)}
=u_e^{\epsilon(e)}\ell_e^{\epsilon(e)}.
\]
The map $\widehat Q_e$ is independent of $x_{|e|}$ and satisfies
$\widehat Q_eu_e^\pm=u_{\bar e}^\pm$. Therefore
\[
\begin{aligned}
[x_{|e|}]B_e
&=\widehat Q_e\,[x_{|e|}]\,C_e(x_{|e|}^{\omega(e)})\\
&=\widehat Q_eP_e^{\epsilon(e)}\\
&=\bigl(\widehat Q_eu_e^{\epsilon(e)}\bigr)
  \ell_e^{\epsilon(e)}\\
&=u_{\bar e}^{\epsilon(e)}\ell_e^{\epsilon(e)}.
\end{aligned}
\]
Since $P_e^{\epsilon(e)}$ has rank one and $\widehat Q_e$
is invertible, this leading coefficient has rank one. Since
\[
[x_{|e|}]B_e
=
u_{\bar e}^{\epsilon(e)}\ell_e^{\epsilon(e)},
\]
the factorization $[x_{|e|}]B_e=i_e\circ p_e$ is obtained by setting
\[
p_e:=\ell_e^{\epsilon(e)},
\qquad
i_e(z):=z\,u_{\bar e}^{\epsilon(e)}.
\]

The edge relations and the rank-one leading coefficients are
now in place. It remains to show that these leading terms
compose nontrivially across every reduced junction. 

\begin{lemma}
\label{lem:transition}
Let $e,f$ be consecutive oriented edges with
\[
t(e)=o(f)=v,
\]
and let $g\in G_v$.  If the segment $egf$ is reduced, then the corresponding scalar
transition coefficient is nonzero.
\end{lemma}

\begin{proof}
We first identify the transition map in Definition \ref{defn-leading-block}(v).
Take the auxiliary space to be the enlarged coefficient field,
regarded as a one-dimensional vector space.  For the junction $egf$, where $t(e)=o(f)=v$, the transition map
is therefore $p_fD_v(g)i_e$. For every scalar $z$, we have
\[
\begin{aligned}
\bigl(p_fD_v(g)i_e\bigr)(z)
&=\ell_f^{\epsilon(f)}
  D_v(g)\bigl(z\,u_{\bar e}^{\epsilon(e)}\bigr)\\
&=z\bigl(\ell_f^{\epsilon(f)}
  D_v(g)u_{\bar e}^{\epsilon(e)}\bigr).
\end{aligned}
\]
Thus the transition map is multiplication by the scalar
\[
\ell_f^{\epsilon(f)}D_v(g)u_{\bar e}^{\epsilon(e)},
\]
and it is invertible precisely when this scalar is nonzero.
We now verify this nonvanishing for every reduced segment $egf$.
Put
\[
d_{f,\bar e}(g)
=
\lambda_f^0\sigma_v(g)Z_{\bar e}^0.
\]
The four possible transition coefficients are obtained directly from the
definitions.  For example,
\[
D_v(g)u_{\bar e}^+
=
b_v+s_{\bar e}\sigma_v(g)Z_{\bar e}^0,
\]
implying
\[
\ell_f^+D_v(g)u_{\bar e}^+
=
\frac{
q_{|f|}
-(s_{\bar e}/s_f)d_{f,\bar e}(g)
}{
q_{|f|}-1
},
\]
and
\[
\ell_f^-D_v(g)u_{\bar e}^+
=
\frac{
-1
+(s_{\bar e}/s_f)d_{f,\bar e}(g)
}{
q_{|f|}-1
}.
\]
The other two expressions are obtained by replacing
$u_{\bar e}^+$ with $u_{\bar e}^-$ and have the same form, with an
additional factor $q_{|e|}$ multiplying $d_{f,\bar e}(g)$.  Explicitly,
they are
\[
\ell_f^+D_v(g)u_{\bar e}^-
=
\frac{
q_{|f|}
-q_{|e|}(s_{\bar e}/s_f)d_{f,\bar e}(g)
}{
q_{|f|}-1
}
\]
and
\[
\ell_f^-D_v(g)u_{\bar e}^-
=
\frac{
-1
+q_{|e|}(s_{\bar e}/s_f)d_{f,\bar e}(g)
}{
q_{|f|}-1
}.
\]

Suppose first that $f=\bar e$.    Hence
\[
s_{\bar e}=s_f,
\qquad
q_{|e|}=q_{|f|}=:q,
\qquad
d_{f,\bar e}(g)=d_{f,f}(g).
\]
Since the segment $egf$ is reduced,
\[
g\notin H_{\bar e}=H_f.
\]
The pointed detector property at the incidence $f$ (see (\ref{eq:pi-dec-prop}) at the beginning of this section) therefore gives
\[
d_{f,f}(g)
=
\lambda_f^0\sigma_v(g)Z_f^0
\neq1.
\]

Moreover, $\epsilon(f)$ and $\epsilon(e)$ have opposite signs.  If
$\epsilon(e)=+$ and $\epsilon(f)=-$, then the transition coefficient is
\[
\ell_f^-D_v(g)u_f^+
=
\frac{d_{f,f}(g)-1}{q-1},
\]
which is nonzero.  If $\epsilon(e)=-$ and $\epsilon(f)=+$, then the
transition coefficient is
\[
\ell_f^+D_v(g)u_f^-
=
\frac{q-q\,d_{f,f}(g)}{q-1}
=
\frac{q\bigl(1-d_{f,f}(g)\bigr)}{q-1},
\]
which is also nonzero.

Suppose instead that $f\neq\bar e$.  Then the  incidences
$\bar e$ and  $f$ are distinct.  This remains
true when $e=f$ is a loop edge, since $e$ and $\bar e$ are distinct
incidences of the loop.  The variables $s_{\bar e}$ and $s_f$ are
therefore algebraically independent.

After clearing the nonzero denominator $q_{|f|}-1$, the four possible
numerators are
\[
q_{|f|}
-\frac{s_{\bar e}}{s_f}d_{f,\bar e}(g),
\qquad
-1+\frac{s_{\bar e}}{s_f}d_{f,\bar e}(g),
\]
and
\[
q_{|f|}
-q_{|e|}\frac{s_{\bar e}}{s_f}d_{f,\bar e}(g),
\qquad
-1
+q_{|e|}\frac{s_{\bar e}}{s_f}d_{f,\bar e}(g).
\]
If $d_{f,\bar e}(g)=0$, these numerators reduce to either
$q_{|f|}$ or $-1$, and hence are nonzero.  If
$d_{f,\bar e}(g)\neq0$, each numerator is a sum of two distinct Laurent
monomials in the independent incidence variables $s_{\bar e}$ and
$s_f$.  Hence none of them can vanish.

Thus every reduced transition coefficient is nonzero.
\end{proof}

We can now obtain the rank-one version of the combination theorem.

\begin{theorem}
    [rank-one combination theorem]
\label{prop:combination}
Let $k$ be a field and let $\cG$ be a finite connected graph of groups with
injective edge maps.  For every vertex $v$, suppose that there is a faithful
representation
\[
 \sigma_v:G_v\longrightarrow\GL(V_v),
 \qquad \dim_{k}V_v=n.
\]
For every incidence $a$ at $v$, suppose that there are
\[
 Z_a^0\in V_v,
 \qquad
 \lambda_a^0\in V_v^*=\Hom_k(V_v,k)
\]
such that:
\begin{enumerate}[label=\textup{(\roman*)}]
\item $\lambda_a^0(Z_a^0)=1$;
\item $H_a$ fixes $Z_a^0$ and preserves $\ker\lambda_a^0$;
\item (pointed detector property) for every $g\in G_v$,
\[
 \lambda_a^0\bigl(\sigma_v(g)Z_a^0\bigr)=1
 \quad\Longleftrightarrow\quad g\in H_a;
\]
\item for each preferred orientation $e$ of a geometric edge, there is a
$k$-linear isomorphism
\[
 Q_e^0:V_{o(e)}\longrightarrow V_{t(e)}
\]
such that, for every $h\in G_{\abs e}$,
\[
 Q_e^0\sigma_{o(e)}(\iota_e(h))
 =\sigma_{t(e)}(\iota_{\bar e}(h))Q_e^0,
\]
and
\[
 Q_e^0Z_e^0=Z_{\bar e}^0,
 \qquad
 \lambda_{\bar e}^0Q_e^0=\lambda_e^0.
\]
\end{enumerate}
Then there is a purely transcendental extension $\K/k$ and an injective
homomorphism
\[
 \rho:\pi_1(\cG,v_*,\cT)\hookrightarrow\GL_{n+1}(\K).
\]

\end{theorem}
\begin{proof}
Let $\mathcal I$ and $\mathcal E$ denote the sets of incidences
and geometric edges, respectively. Put
\[
F=k(s_a:a\in\mathcal I,\ q_{|e|}:|e|\in\mathcal E),
\qquad
K=F(x_{|e|}:|e|\in\mathcal E),
\]
where all the displayed indeterminates are algebraically
independent over $k$.

After extending the enlarged vertex spaces and representations
to $F$, the preceding construction gives a leading-block system
of rank one and vertex dimension $n+1$ over $F$, with Laurent
variables $x_{|e|}$. Lemma~2.2 verifies the transition condition.
Applying Theorem \ref{thm:leading-block} with coefficient field $F$ therefore gives
\[
\pi_1(\mathcal G,v_*,T)\hookrightarrow\mathrm{GL}_{n+1}(K).
\]
Since $K/k$ is purely transcendental, the conclusion follows.
\end{proof}

We now discuss two special cases of Theorem \ref{prop:combination} for amalgamated products and HNN extensions.

\begin{corollary}
    [Amalgamated criterion]
\label{cor:amalgamated-matrix-entry}
Let $k$ be a field, let $V$ be an $n$-dimensional
$\Bbbk$-vector space, and let
\[
\sigma_A\colon A\longrightarrow \GL(V),
\qquad
\sigma_B\colon B\longrightarrow \GL(V)
\]
be faithful representations of two groups $A,B$. Suppose that $H$ is a common subgroup of
$A$ and $B$ and that
\[
\sigma_A|_H=\sigma_B|_H.
\]
Assume that there exist $Z\in V$ and $\lambda\in V^*$ such that
\[
\lambda(Z)=1,
\]
the subgroup $H$ fixes the pair $(Z,\lambda)$, in the sense that
\[
\sigma_A(h)Z=Z,
\qquad
\lambda\sigma_A(h)=\lambda
\quad(h\in H),
\]
and, for each $X\in\{A,B\}$,
\[
\lambda\sigma_X(g)Z=1
\quad\Longleftrightarrow\quad
g\in H.
\]
Then there is a purely transcendental extension
$\mathbb F/k$ for which
\[
A*_H B\hookrightarrow \GL_{n+1}(\mathbb F).
\]
In particular, $A*_H B$ is linear.
\end{corollary}

\begin{proof}
Regard $A*_H B$ as the fundamental group of the graph of groups
$\mathcal G$ having two vertices, with vertex groups $A$ and $B$, and
one geometric edge, with edge group $H$.

Assign $\sigma_A$ and $\sigma_B$ to the two vertex groups. Choose an
orientation $e$ from the $A$-vertex to the $B$-vertex. At the two
incidences of the unique edge, put
\[
Z_e^0=Z_{\bar e}^0=Z,
\qquad
\lambda_e^0=\lambda_{\bar e}^0=\lambda,
\]
and define
\[
Q_e^0=\id_V.
\]
Since
\[
\sigma_A|_H=\sigma_B|_H,
\]
for every $h\in H$ we have
\[
Q_e^0\sigma_A(h)
=
\sigma_B(h)Q_e^0.
\]
Moreover,
\[
Q_e^0Z_e^0=Z_{\bar e}^0,
\qquad
\lambda_{\bar e}^0Q_e^0=\lambda_e^0.
\]
Thus condition \textup{(iv)} of Theorem~\ref{prop:combination} is
satisfied.

For every $h\in H$,
\[
\sigma_A(h)Z=\sigma_B(h)Z=Z,
\qquad
\lambda\sigma_A(h)=\lambda\sigma_B(h)=\lambda.
\]
Thus the edge group fixes the distinguished vector and covector at both
ends of the edge.

Finally, the hypothesis
\[
\lambda\sigma_X(g)Z=1
\quad\Longleftrightarrow\quad g\in H
\qquad(X\in\{A,B\})
\]
is exactly the pointed detector property in
Theorem~\ref{prop:combination}. All the hypotheses of that theorem are
therefore satisfied. It supplies a faithful representation
\[
\pi_1(\mathcal G)\longrightarrow \GL_{n+1}(\mathbb F)
\]
over the rational function field $\mathbb F$ generated by the auxiliary
parameters in its construction. Since
\[
\pi_1(\mathcal G)\cong A*_H B,
\]
the required embedding follows.
\end{proof}

\begin{remark}
In the preceding corollary, choose a basis of $V$ consisting of $Z$
followed by a basis of $\ker\lambda$. Then $\lambda$ is the first
coordinate functional. The requirement that $H$ fix both $Z$ and
$\lambda$ is equivalent to
\[
\sigma_X(H)\subseteq \{1\}\times\mathrm{GL}_{n-1}(k)
\qquad (X\in\{A,B\}),
\]
where the product is embedded block-diagonally in $\mathrm{GL}_n(k$).
Moreover,
\[
\lambda\sigma_X(g)Z=(\sigma_X(g))_{11},
\]
so the pointed detector property becomes
\[
g\in H
\quad\Longleftrightarrow\quad
(\sigma_X(g))_{11}=1.
\]
Thus the criterion requires that the two faithful representations
agree on $H$, that their common restriction to $H$ have the above
block-diagonal form, and that membership in $H$ be detected by the
top-left matrix entry. This recovers the linearity conclusion of
Douba--Tsouvalas~\cite[Theorem~2.1(1)]{DoubaTs}.

\end{remark}

\begin{corollary}[HNN criterion]
\label{cor:hnn-matrix-entry}
Let $k$ be a field, let $V$ be an $n$-dimensional
$\Bbbk$-vector space, and let
\[
\sigma\colon G\longrightarrow \GL(V)
\]
be faithful. Let $H,K\leq G$, and let
$\phi\colon H\xrightarrow{\sim}K$ be an isomorphism. Suppose that there are nonzero vectors
\[
Z_H,Z_K\in V,
\]
linear functionals
\[
\lambda_H,\lambda_K\in V^*,
\]
and a matrix $Q\in\GL(V)$ such that
\[
\lambda_H(Z_H)=\lambda_K(Z_K)=1,
\]
and the following conditions hold:
\begin{enumerate}
\item the associated subgroups fix the corresponding vectors and
covectors:
\[
\sigma(h)Z_H=Z_H,\qquad
\lambda_H\sigma(h)=\lambda_H
\quad(h\in H),
\]
and
\[
\sigma(k)Z_K=Z_K,\qquad
\lambda_K\sigma(k)=\lambda_K
\quad(k\in K);
\]

\item each associated subgroup is detected by its corresponding matrix
coefficient:
\[
\lambda_H\sigma(g)Z_H=1
   \quad\Longleftrightarrow\quad g\in H,
\qquad
\lambda_K\sigma(g)Z_K=1
   \quad\Longleftrightarrow\quad g\in K;
\]

\item $Q$ realizes the attaching isomorphism and identifies the pointed
data:
\[
Q\sigma(h)Q^{-1}=\sigma(\phi(h)),\qquad
QZ_H=Z_K,\qquad
\lambda_KQ=\lambda_H
\quad(h\in H).
\]
\end{enumerate}
Then there is a purely transcendental extension
$\mathbb F/k$ such that
\[
G*_{\phi}
 =
\left\langle
G,u\ \middle|\ uhu^{-1}=\phi(h)\ \text{for all }h\in H
\right\rangle
\hookrightarrow \GL_{n+1}(\mathbb F).
\]
In particular, the HNN extension $G*_{\phi}$ is linear.
\end{corollary}

\begin{proof}
Regard $G*_{\phi}$ as the fundamental group of the graph of groups
$\mathcal G$ having one vertex, with vertex group $G$, and one geometric
loop edge, with edge group $H$. Choose an orientation $e$ of the loop
and define the two edge monomorphisms by
\[
\iota_e(h)=h,
\qquad
\iota_{\bar e}(h)=\phi(h).
\]
Thus the two incidences $e$ and $\bar e$ correspond respectively to the
subgroups
\[
H_e=H,
\qquad
H_{\bar e}=K.
\]

Assign the faithful representation $\sigma$ to the unique vertex group.
At the two incidences put
\[
Z_e^0=Z_H,
\qquad
\lambda_e^0=\lambda_H,
\qquad
Z_{\bar e}^0=Z_K,
\qquad
\lambda_{\bar e}^0=\lambda_K.
\]
The normalization assumptions verify condition~\textup{(i)} of
Theorem~\ref{prop:combination}. By condition~\textup{(1)}, each
incident edge subgroup fixes the corresponding distinguished vector and
covector, and hence preserves the kernel of the corresponding
covector. Thus condition~\textup{(ii)} is satisfied. Condition~\textup{(2)}
is exactly the pointed detector property hypothesis in
condition~\textup{(iii)} of Theorem~\ref{prop:combination}.

Finally, put
\[
Q_e^0=Q.
\]
For every $h\in H$, condition~\textup{(3)} gives
\[
\begin{aligned}
Q_e^0\sigma\bigl(\iota_e(h)\bigr)
&=Q\sigma(h)\\
&=\sigma\bigl(\phi(h)\bigr)Q\\
&=\sigma\bigl(\iota_{\bar e}(h)\bigr)Q_e^0.
\end{aligned}
\]
Moreover,
\[
Q_e^0Z_e^0=QZ_H=Z_K=Z_{\bar e}^0,
\qquad
\lambda_{\bar e}^0Q_e^0
=\lambda_KQ
=\lambda_H
=\lambda_e^0.
\]
Thus condition~\textup{(iv)} of
Theorem~\ref{prop:combination} is satisfied.

All the hypotheses of Theorem~\ref{prop:combination} are therefore
satisfied. It yields a faithful representation
\[
\pi_1(\mathcal G)\longrightarrow \GL_{n+1}(\mathbb F),
\]
where $\mathbb F/\Bbbk$ is the purely transcendental extension generated
by the auxiliary parameters in the construction. 
\end{proof}

\section{Product Seifert blocks and the graph-manifold application}
\label{sec:product-blocks}

We now verify the hypotheses of \Cref{prop:combination} for a graph manifold
whose Seifert blocks are trivial circle bundles over orientable surfaces of positive genus. In this case, we have a graph of groups with each vertex group a product of a surface group and $\mathbb{Z}$, and each edge group isomorphic to $\mathbb{Z}^2$ (called peripheral subgroup).  The construction is organized so that the
two generators of each peripheral subgroup can be rescaled independently
enough to accommodate an arbitrary torus gluing.

So let us fix such a graph of groups from now on, and pick a preferred orientation for the underlying graph, i.e., pick an oriented edge for each geometric edge.

\subsection{Representations of the vertex groups}\label{Section:rep-prod-group} We first construct  representations of the vertex groups. Let
\[
M_v=Q_v\times S^1,
\]
where $Q_v$ is a compact connected orientable surface of positive genus with
nonempty boundary.  Let $t_v$ denote the positive generator of the circle
factor.  If $Q_v$ has genus $h\geq1$ and boundary components represented by
$c_1,\dots,c_b$, choose the presentation
\begin{equation}
\label{eq:surface-presentation}
\pi_1(Q_v)=
\left\langle
 a_1,b_1,\dots,a_h,b_h,c_1,\dots,c_b
 \ \middle|\ 
 \prod_{j=1}^h[a_j,b_j]c_1\cdots c_b=1
\right\rangle.
\end{equation}
Then
\[
G_v:=\pi_1(M_v)=\pi_1(Q_v)\times\langle t_v\rangle,
\]
and the peripheral subgroup associated with the $i$th boundary component is
\[
H_{v,i}=\langle c_i,t_v\rangle\cong\Z^2.
\]

\subsubsection*{Fuchsian representation}

Choose a complete finite-area hyperbolic metric on the interior of $Q_v$ and a
faithful holonomy representation
\[
\bar\rho_v\colon\pi_1(Q_v)\longrightarrow\PSL_2(\R)
\]
for which each $c_i$ is represented by a primitive parabolic.  Such a metric
and holonomy are standard; see, for example, \cite[Chapters~7--9]{Beardon}.
Since $\pi_1(Q_v)$ is free, $\bar\rho_v$ lifts to
\[
\rho_v\colon\pi_1(Q_v)\longrightarrow\SL_2(\R).
\]
Write
\begin{equation}
\label{eq:boundary-lifts}
\rho_v(c_i)=\epsilon_iU_i,
\qquad \epsilon_i\in\{\pm1\},
\end{equation}
where $U_i\neq I_2$ is unipotent with trace $2$. Note that each $U_i$ is individually conjugate to an upper triangular
unipotent matrix.

We use a second two-dimensional representation to cancel the signs in
\eqref{eq:boundary-lifts}.  There is a homomorphism
\[
\eta_v\colon\pi_1(Q_v)\longrightarrow\GL_2(\R)
\]
with $\eta_v(c_i)=\epsilon_iI_2$.  Indeed, if
$\prod_i\epsilon_i=1$, send all handle generators to $I_2$.  If
$\prod_i\epsilon_i=-1$, put
\[
R_0=\begin{pmatrix}1&0\\0&-1\end{pmatrix},
\qquad
J_0=\begin{pmatrix}0&1\\1&0\end{pmatrix}.
\]
Then $[R_0,J_0]=-I_2$, so sending $a_1\mapsto R_0$, $b_1\mapsto J_0$, and all
remaining handle generators to $I_2$ makes relation
\eqref{eq:surface-presentation} hold.

We record the following elementary fact on tensor product which will be used below.

\begin{lemma}[Scalar tensor products]
\label{lem:scalar-tensor}
Let $A\in\GL_r(k)$ and $B\in\GL_s(k)$, where $r,s\geq1$.  If $A\otimes B$ is
a scalar matrix, then both $A$ and $B$ are scalar matrices.
\end{lemma}

\begin{proof}
Write $A=(a_{ij})$ and view $A\otimes B$ as an $r\times r$ block matrix whose
$(i,j)$-block is $a_{ij}B$.  Since $B$ is invertible, every off-diagonal block
vanishes only if $a_{ij}=0$ for $i\neq j$.  If
$A\otimes B=\lambda I_{rs}$, each diagonal block satisfies
$a_{ii}B=\lambda I_s$.  Thus every $a_{ii}$ is nonzero, $B$ is scalar, and
comparison of the diagonal blocks shows that all $a_{ii}$ are equal.  Hence
$A$ is scalar as well.
\end{proof}

Define a homomorphism $\alpha_v: \pi_1(Q_v)\times \langle t_v\rangle \to \GL_6(\mathbb{R})$ by 
\begin{equation}
\label{eq:alpha}
\alpha_v(wt_v^m)
 :=\bigl(\rho_v(w)\otimes\eta_v(w)\bigr)\oplus I_2
 \in\GL_6(\R).
\end{equation}
In particular, $\alpha_v(t_v)=I_6$.  The restriction of $\alpha_v$ to
$\pi_1(Q_v)$ is faithful.  Indeed, if $\alpha_v(w)=I_6$, then
$\rho_v(w)\otimes\eta_v(w)=I_4$.  By \Cref{lem:scalar-tensor}, $\rho_v(w)$ is
scalar; faithfulness of $\bar\rho_v$ then gives $w=1$.

The sign cancellation gives the exact peripheral matrix
\begin{equation}
\label{eq:Kvi}
K_{v,i}:=\alpha_v(c_i)=(U_i\otimes I_2)\oplus I_2.
\end{equation}
Thus $K_{v,i}-I_6$ is nonzero, has square zero, and has rank $2$.

\subsubsection*{Elementary matrices with adjustable slopes} 

Identify $\Mat_6(\R)$ with $\Mat_3(\R)\otimes\Mat_2(\R)$, and let $E_{ij}$
denote the $3\times3$ matrix whose $(i,j)$-th entry is $1$ and all other
entries are zero.  Put
\begin{equation}
\label{eq:elementary-matrices}
K=I_6+E_{13}\otimes I_2,
\qquad
X=I_6+E_{12}\otimes I_2,
\qquad
Y(r)=I_6+rE_{23}\otimes I_2.
\end{equation}
For $r\in\R$, write
\[
K^r:=I_6+r(K-I_6).
\]
Since $(K-I_6)^2=0$, these matrices form an additive one-parameter subgroup:
\[
K^rK^s=K^{r+s},
\qquad
(K^r)^{-1}=K^{-r}.
\]
For integral $r$ this agrees with the ordinary matrix power.  Direct
multiplication gives
\begin{equation}
\label{eq:Heisenberg-commutators}
[X,Y(r)]=K^r,
\qquad
[X,K]=[Y(r),K]=I_6.
\end{equation}

Assign a rational number $\mu_{v,i}$ to each boundary curve $c_i$ of $Q_v$.  The
values $\mu_{v,i}$ belonging to glued JSJ tori will be specified from the gluing matrices
in \Cref{sec:arbitrary-gluings}. Put
\[
\mu_v:=\sum_{i=1}^b\mu_{v,i}.
\]
Define
\[
\beta_v\colon G_v\longrightarrow\GL_6(\R)
\]
by
\begin{align}
\beta_v(c_i)&=K^{\mu_{v,i}},
&\beta_v(t_v)&=K,
\\
\beta_v(a_1)&=X,
&\beta_v(b_1)&=Y(-\mu_v),
\label{eq:beta-handle}
\end{align}
and by sending all remaining handle generators to $I_6$.  Equations
\eqref{eq:Heisenberg-commutators} give
\[
[X,Y(-\mu_v)]\prod_{i=1}^bK^{\mu_{v,i}}
 =K^{-\mu_v}K^{\mu_v}=I_6,
\]
so the surface relation is satisfied.  The image of $\pi_1(Q_v)$ centralizes
$K$, and hence these assignments extend to a homomorphism on the direct
product $G_v$.

\subsection{Tensor scaling}
\label{subsec:pointed-scaling} In this subsection, we further modify the representations of the vertex groups to provide the flexibility needed for the subsequent gluing construction.

Let $k$ be a field of characteristic zero, and regard $\Mat_d(k)$ as a
$d^2$-dimensional $k$-vector space.  We distinguish the multiplicative
conjugation operator
\[
\Ad(B)(Z):=BZB^{-1}
\qquad(B\in\GL_d(k))
\]
from the additive commutator operator
\[
\ad(X)(Z):=XZ-ZX
\qquad(X\in\Mat_d(k)).
\]

Let
\[
A=I_d+N\in\GL_d(k),
\]
for some $N$ satisfying $\qquad N^2=0$. For $r\in k$ and $Z\in\Mat_d(k)$, define
\begin{equation}
\label{eq:UA-tau}
A^r:=I_d+rN,
\qquad
U_A(r):=\Ad(A^r),
\qquad
\tau_A(Z):=\frac1d\tr(ZA^{-1}).
\end{equation}
For arbitrary $r\in k$, the notation $A^r$ is a definition; for
$r\in\Z$ it agrees with the ordinary matrix power.

Put $D:=\ad(N)$.  For every $Z\in\Mat_d(k)$,
\begin{align}
D^2(Z)&=N(NZ-ZN)-(NZ-ZN)N=-2NZN,
\label{eq:D-square}\\
D^3(Z)&=0.
\label{eq:D-cube}
\end{align}
Thus all exponentials below are finite polynomials, and
\begin{align*}
\exp(rD)(Z)
&=(\sum_{n\geq 0} \frac{(rD)^n}{n!})(Z) \\
 &=Z+r(NZ-ZN)-r^2NZN\\
 &=(I_d+rN)Z(I_d-rN)\\
 &=A^rZ(A^r)^{-1}.
\end{align*}
Consequently,
\begin{equation}
\label{eq:U-exponential}
U_A(r)=\exp(rD).
\end{equation}

\begin{lemma}[Trace scaling]
\label{lem:pointed-scaling}
Let $A=I_d+N\in\GL_d(k)$ with $N^2=0$, and let $c\in k^\times$.  There is a
$k$-linear automorphism $S_c\in\GL(\Mat_d(k))$ such that, for every $r\in k$,
\begin{equation}
\label{eq:pointed-scaling}
S_c(A)=A,
\qquad
\tau_A\circ S_c=\tau_A,
\qquad
S_cU_A(r)S_c^{-1}=U_A(cr).
\end{equation}

\end{lemma}

\begin{remark}
    The automorphism $S_c$ is noncanonical and depends on $c$; only its existence is used later.
\end{remark}

\begin{proof}
Note that the operator $D=\ad(N)$ is nilpotent.

\smallskip
\noindent\emph{Step 1: a $D$-invariant decomposition.}
Because $A$ is a polynomial in $N$, one has $D(A)=\ad(N)(A)=0$.  Also $A^{-1}=I_d-N$
commutes with $N$, and the commutativity of trace gives, for every $Z\in \Mat_d(k)$, that
\begin{align*}
\,\tau_A(DZ)
 &=\frac{1}{d}\tr\bigl((NZ-ZN)A^{-1}\bigr)\\
 &=\frac{1}{d}(\tr(ZA^{-1}N)-\tr(ZNA^{-1}))=0.
\end{align*}
Hence $\tau_A\circ D=0$.  Moreover $\tau_A(A)=1$, so we have a decomposition
\begin{equation}
\label{eq:pointed-decomposition}
\Mat_d(k)=kA\oplus W,
\qquad W:=\ker\tau_A.
\end{equation}
Both summands are $D$-invariant. Indeed, $D(A)=0$, so $D$
vanishes on $kA$, while $\tau_A\circ D=0$ implies that
$D(W)\subseteq W$.

\smallskip
\noindent\emph{Step 2: rescaling the Jordan chains of $D|_W$.}
By \eqref{eq:D-cube}, the nilpotent operator $D|_W$ satisfies
$D^3=0$.  Its minimal polynomial divides $x^3$, so it splits over $k$, and
$W$ has a basis partitioned into nilpotent Jordan chains of length at most
three.  On a chain
\[
v_0,v_1,\dots,v_{\ell-1},
\qquad
Dv_0=0,
\quad Dv_j=v_{j-1}\ (1\leq j<\ell),
\]
define
\[
S_c(v_j):=c^{-j}v_j.
\]
Then $S_cD|_WS_c^{-1}=cD|_W$ on that chain.  Doing this on every chain gives an
automorphism of $W$ with the same property.  Extend it to all of
$\Mat_d(k)$ by making it the identity on $kA$.  Then
\begin{equation}
\label{eq:S-conjugates-D}
S_cDS_c^{-1}=cD,
\qquad
S_c(A)=A.
\end{equation}

\smallskip
\noindent\emph{Step 3: preservation of the trace functional.}
Write $Z=\lambda A+w$ according to \eqref{eq:pointed-decomposition}.  Since
$S_c$ fixes $A$ and preserves $W$,
\[
\tau_A(S_cZ)=\tau_A(\lambda A+S_cw)=\lambda=\tau_A(Z).
\]

\smallskip
\noindent\emph{Step 4: scaling the conjugation parameter.}
Using \eqref{eq:U-exponential}, \eqref{eq:S-conjugates-D}, and the fact that
the exponential is a finite polynomial in $D$,
\[
S_cU_A(r)S_c^{-1}
 =S_c\exp(rD)S_c^{-1}
 =\exp(rcD)=U_A(cr).
\]
This proves \eqref{eq:pointed-scaling}.
\end{proof}

We apply the lemma to the standard unipotent matrix $K$ from
\eqref{eq:elementary-matrices}.  Write
\begin{equation}
\label{eq:U-tau-K}
U(r):=\Ad(K^r),
\qquad
\tau_K(Z):=\frac16\tr(ZK^{-1}).
\end{equation}
The matrices $K_{v,i}-I_6$ defined by \eqref{eq:Kvi} and $K-I_6$ defined by \eqref{eq:elementary-matrices} have the same Jordan form: both square
to zero and have rank $2$, hence both have Jordan type $(2,2,1,1)$.  Choose
$T_{v,i}\in\GL_6(\R)$ such that
\[
T_{v,i}K_{v,i}T_{v,i}^{-1}=K,
\]
and set
\[
L_{v,i}:=\Ad(T_{v,i}).
\]
Then
\begin{equation}
\label{eq:standardization}
L_{v,i}(K_{v,i})=K,
\qquad
\tau_K\circ L_{v,i}=\tau_{K_{v,i}}.
\end{equation}
Indeed, for $Z\in\Mat_6(\R)$,
\begin{align*}
\tau_K(L_{v,i}(Z))
 &=\frac16\tr\bigl(T_{v,i}ZT_{v,i}^{-1}K^{-1}\bigr)\\
 &=\frac16\tr\bigl(ZT_{v,i}^{-1}K^{-1}T_{v,i}\bigr)\\
 &=\frac16\tr(ZK_{v,i}^{-1}).
\end{align*}

On the $72$-dimensional vector space
\[
V_v:=\Mat_6(\R)\oplus\Mat_6(\R)
\]
define
\begin{equation}
\label{eq:Theta}
\sigma_v:=(\Ad\circ\alpha_v)\oplus(\Ad\circ\beta_v)
\colon G_v\longrightarrow\GL(V_v),
\end{equation}
where $\alpha_v$ is defined in \eqref{eq:alpha}  and $\beta_v$ is defined in \eqref{eq:beta-handle}. 

\begin{lemma}
\label{lem:Theta-faithful}
The representation $\sigma_v$ is faithful.
\end{lemma}

\begin{proof}
Suppose that $g=wt_v^m$ lies in its kernel.  Then
$\Ad\alpha_v(w)=I$.  A matrix whose conjugation action on the full matrix
algebra is trivial commutes with every matrix and is therefore scalar; hence
$\alpha_v(w)$ is scalar.  The final trivial $I_2$-summand in
\eqref{eq:alpha} forces that scalar to be $1$.  Faithfulness of
$\alpha_v|_{\pi_1(Q_v)}$ gives $w=1$.

The second block now says $\Ad(K^m)=I$, so the same centralizer argument makes
$K^m$ scalar.  But
\[
K^m=I_6+m(K-I_6),
\]
and $K-I_6$ is a nonzero off-diagonal nilpotent matrix.  Over characteristic
zero this is scalar only when $m=0$.
\end{proof}

\subsection{Torus gluings}
\label{sec:arbitrary-gluings} We now start to ``glue'' the representations together, i.e., we define the maps $Q_e^0$ in Theorem \ref{prop:combination} (iv) for each preferred oriented edge we choose.
Fix one such oriented edge $e$ and put
\[
v_-:=o(e),\qquad v_+:=t(e).
\]
Let $a_-$ be the incidence corresponding to $e$ and $a_+$ to $\bar{e}$.  Write
\[
H_{a_-}=\langle c_-,t_-\rangle\leq G_{v_-},
\qquad
H_{a_+}=\langle c_+,t_+\rangle\leq G_{v_+},
\]
where \(c_\pm\) are the elements represented by the corresponding boundary curves of \(Q_{v_\pm}\), $t_-=t_{v_-}$ and $t_+=t_{v_+}$ are the Seifert-fiber generators in the
two product blocks.

Choose the source peripheral basis as a basis of $G_{|e|}\cong\Z^2$.  If
$h_{p,q}\in G_{|e|}$ corresponds to the exponent column $(p,q)^T \in \Z^2$, define the two
edge inclusions by
\begin{equation}
\label{eq:edge-inclusions}
\iota_e(h_{p,q})=c_-^pt_-^q,
\qquad
\iota_{\bar e}(h_{p,q})=c_+^{ap+bq}t_+^{cp+dq},
\end{equation}
where
\begin{equation}
\label{eq:gluing-matrix}
A=\begin{pmatrix}a&b\\c&d\end{pmatrix}\in\GL_2(\Z),
\qquad
\begin{pmatrix}p\\q\end{pmatrix}
\longmapsto A\begin{pmatrix}p\\q\end{pmatrix},
\qquad
\Delta:=\det A\in\{\pm1\}.
\end{equation}

For $\epsilon\in\{-,+\}$, let $i_\epsilon$ be the index of the
boundary component represented by $c_\epsilon$ in $Q_{v_\epsilon}$.
Put
\[
K_\epsilon:=K_{v_\epsilon,i_\epsilon}
           =\alpha_{v_\epsilon}(c_\epsilon),
\qquad
T_\epsilon:=T_{v_\epsilon,i_\epsilon}.
\]
By the choice of $T_{v_\epsilon,i_\epsilon}$,
\[
T_\epsilon K_\epsilon T_\epsilon^{-1}=K.
\]
Define the linear maps $L_\epsilon$ and $J_\epsilon$ by
\[
L_\epsilon(P):=T_\epsilon PT_\epsilon^{-1},
\qquad
J_\epsilon(P,Q):=(T_\epsilon PT_\epsilon^{-1},Q),
\qquad P,Q\in\Mat_6(\R).
\]
Thus $L_\epsilon=\Ad(T_\epsilon)$ and
$J_\epsilon=L_\epsilon\oplus\id_{\Mat_6(\R)}$.

For every integer $p$, write
\[
L_\epsilon\Ad(K_\epsilon^p)L_\epsilon^{-1}
 =\Ad(K^p)=:U(p).
\]
Recall that
\[
\sigma_{v_\epsilon}
=
(\Ad\circ\alpha_{v_\epsilon})
\oplus
(\Ad\circ\beta_{v_\epsilon}),
\]
and that the peripheral generators satisfy
\[
\alpha_{v_\epsilon}(c_\epsilon)=K_\epsilon,
\qquad
\alpha_{v_\epsilon}(t_\epsilon)=I_6,
\qquad
\beta_{v_\epsilon}(c_\epsilon)=K^{\mu_\epsilon},
\qquad
\beta_{v_\epsilon}(t_\epsilon)=K.
\]
Since $\alpha_{v_\epsilon}$ and $\beta_{v_\epsilon}$ are homomorphisms,
for every $p,q\in\mathbb Z$ we have
\[
\alpha_{v_\epsilon}(c_\epsilon^pt_\epsilon^q)
=
K_\epsilon^p,
\qquad
\beta_{v_\epsilon}(c_\epsilon^pt_\epsilon^q)
=
(K^{\mu_\epsilon})^pK^q
=
K^{\mu_\epsilon p+q}.
\]
Consequently,
\[
\sigma_{v_\epsilon}(c_\epsilon^pt_\epsilon^q)
=
\Ad(K_\epsilon^p)\oplus\Ad(K^{\mu_\epsilon p+q}).
\]
So for every $P\in\Mat_6(\mathbb R)$,
\[
\begin{aligned}
\bigl(L_\epsilon\Ad(K_\epsilon^p)L_\epsilon^{-1}\bigr)(P)
&=
T_\epsilon
\bigl(
K_\epsilon^p(T_\epsilon^{-1}PT_\epsilon)K_\epsilon^{-p}
\bigr)
T_\epsilon^{-1}\\
&=
(T_\epsilon K_\epsilon^pT_\epsilon^{-1})
P
(T_\epsilon K_\epsilon^{-p}T_\epsilon^{-1})\\
&=
K^pPK^{-p}\\
&=
\Ad(K^p)(P).
\end{aligned}
\]
Thus
\[
L_\epsilon\Ad(K_\epsilon^p)L_\epsilon^{-1}
=
\Ad(K^p)
=
U(p).
\]
Now define
\begin{equation}
\label{eq:standardized-peripheral-action}
\begin{aligned}
\widehat\sigma_\epsilon(p,q)
:=
J_\epsilon
\sigma_{v_\epsilon}(c_\epsilon^pt_\epsilon^q)
J_\epsilon^{-1}.
\end{aligned}
\end{equation}
Since $J_\epsilon$ acts as the identity on the second matrix summand,
we obtain
\begin{equation}
\label{eq:standardized-peripheral-action-2}
\begin{aligned}
\widehat\sigma_\epsilon(p,q)
&=
\bigl(
L_\epsilon\Ad(K_\epsilon^p)L_\epsilon^{-1}
\bigr)
\oplus
\Ad(K^{\mu_\epsilon p+q})\\
&=
U(p)\oplus U(\mu_\epsilon p+q).
\end{aligned}
\end{equation}

To formulate the pointed part of the edge condition, set
\begin{equation*}
\label{eq:standard-pointed-data}
\widehat Z:=(K,K),
\qquad
\widehat\lambda(P,Q)
 :=\frac12\bigl(\tau_K(P)+\tau_K(Q)\bigr)
 =\frac{\tr(PK^{-1})+\tr(QK^{-1})}{12}.
\end{equation*}
At the incidence $a_\epsilon$, put
\begin{equation}
\label{eq:unstandardized-pointed-data}
Z_{a_\epsilon}^0:=J_\epsilon^{-1}\widehat Z=(K_\epsilon,K),
\qquad
\lambda_{a_\epsilon}^0:=\widehat\lambda\circ J_\epsilon,
\end{equation}
Thus $\lambda_{a_\epsilon}^0(P,Q)= \frac{\tr(PK_\epsilon^{-1})+\tr(QK^{-1})}{12}.$

\qquad

\begin{lemma}[Arbitrary torus gluing]
\label{lem:arbitrary-gluing}
For the gluing matrix $A$ in \eqref{eq:gluing-matrix}, there are rational
numbers $\mu_-,\mu_+$ and an $\R$-linear isomorphism
\[
Q_e^0\colon V_{v_-}\longrightarrow V_{v_+}
\]
such that, for every $h\in G_{|e|}$,
\begin{equation}
\label{eq:gluing-intertwining}
Q_e^0\sigma_{v_-}(\iota_e(h))
 =\sigma_{v_+}(\iota_{\bar e}(h))Q_e^0,
\end{equation}
and
\begin{equation}
\label{eq:gluing-pointed}
Q_e^0Z_{a_-}^0=Z_{a_+}^0,
\qquad
\lambda_{a_+}^0Q_e^0=\lambda_{a_-}^0.
\end{equation}
More precisely:
\begin{enumerate}[label=\textup{(\alph*)}]
\item if $b\neq0$, take
\begin{equation}
\label{eq:mu-b-nonzero}
\mu_-:=\frac ab,
\qquad
\mu_+:=-\frac db,
\end{equation}
and in standardized coordinates define
\begin{equation}
\label{eq:Qhat-b-nonzero}
\widehat Q_e^0(P,Q):=
\bigl(S_bQ,S_{-\Delta/b}P\bigr);
\end{equation}
\item if $b=0$, take
\begin{equation}
\label{eq:mu-b-zero}
\mu_-:=\frac cd,
\qquad
\mu_+:=0,
\end{equation}
and define
\begin{equation}
\label{eq:Qhat-b-zero}
\widehat Q_e^0(P,Q):=
\bigl(S_aP,S_dQ\bigr).
\end{equation}
\end{enumerate}
Here $S_r$ is any automorphism supplied by \Cref{lem:pointed-scaling}, applied
to $K$, for the indicated nonzero scalar $r$.  In both cases one may take
\begin{equation}
\label{eq:Q-unstandardize}
Q_e^0:=J_+^{-1}\widehat Q_e^0J_-.
\end{equation}
\end{lemma}

\begin{proof}
We first work in the standardized coordinates introduced in
\eqref{eq:standardized-peripheral-action}. For $p,q\in\mathbb Z$, write
\[
\widehat\sigma_-(p,q)
=
U(p)\oplus U(\mu_-p+q)
\]
and define the standardized target action in source coordinates by
\[
\widehat\sigma_+^A(p,q)
:=
\widehat\sigma_+(ap+bq,cp+dq)
=
U(ap+bq)\oplus
U\bigl(\mu_+(ap+bq)+cp+dq\bigr).
\]
Thus $\widehat\sigma_+^A=\widehat\sigma_+\circ A$; the superscript
$A$ denotes precomposition with the gluing map on $\mathbb Z^2$.

We shall use the following consequence of Lemma~\ref{lem:pointed-scaling}.
For every $r\in\mathbb R^\times$, there is an automorphism
\[
S_r\in\GL(\Mat_6(\mathbb R))
\]
such that, for every $u\in\mathbb R$,
\begin{equation}
\label{eq:scaling-identity}
S_rU(u)=U(ru)S_r,
\qquad
S_rK=K,
\qquad
\tau_KS_r=\tau_K.
\end{equation}

\smallskip
\noindent
\emph{Case 1: $b\neq0$.}
Set
\[
\mu_-=\frac{a}{b},
\qquad
\mu_+=-\frac{d}{b},
\]
and introduce the source parameters
\[
x:=p,
\qquad
y:=\mu_-p+q=\frac{a}{b}p+q.
\]
Then
\[
ap+bq=by,
\]
whereas
\[
\begin{aligned}
\mu_+(ap+bq)+cp+dq
&=
-\frac{d}{b}(ap+bq)+cp+dq\\
&=
\frac{bc-ad}{b}\,p
=
-\frac{\Delta}{b}\,x.
\end{aligned}
\]
Consequently,
\[
\widehat\sigma_-(p,q)=U(x)\oplus U(y)
\]
and
\[
\widehat\sigma_+^A(p,q)
=
U(by)\oplus U\left(-\frac{\Delta}{b}x\right).
\]

Since $b\neq0$ and $\Delta=\pm1$, both $b$ and $-\Delta/b$ are
nonzero. Define
\[
\widehat Q_e^0(P,Q)
:=
\left(S_bQ,S_{-\Delta/b}P\right).
\]
Using \eqref{eq:scaling-identity}, we obtain
\[
\begin{aligned}
\widehat Q_e^0\widehat\sigma_-(p,q)(P,Q)
&=
\widehat Q_e^0\bigl(U(x)P,U(y)Q\bigr)\\
&=
\left(S_bU(y)Q,S_{-\Delta/b}U(x)P\right)\\
&=
\left(
U(by)S_bQ,
U\left(-\frac{\Delta}{b}x\right)S_{-\Delta/b}P
\right)\\
&=
\widehat\sigma_+^A(p,q)\widehat Q_e^0(P,Q).
\end{aligned}
\]
The map $\widehat Q_e^0$ is invertible, since it is the composite of
the block interchange with the invertible map
$S_{-\Delta/b}\oplus S_b$.

Moreover, writing
\[
\widehat Z=(K,K),
\qquad
\widehat\lambda(P,Q)
=
\frac{1}{2}\bigl(\tau_K(P)+\tau_K(Q)\bigr),
\]
the remaining identities in \eqref{eq:scaling-identity} give
\[
\widehat Q_e^0\widehat Z
=
(S_bK,S_{-\Delta/b}K)
=
(K,K)
=
\widehat Z
\]
and
\[
\begin{aligned}
\widehat\lambda\bigl(\widehat Q_e^0(P,Q)\bigr)
&=
\frac{1}{2}
\left(
\tau_K(S_bQ)+\tau_K(S_{-\Delta/b}P)
\right)\\
&=
\frac{1}{2}\bigl(\tau_K(Q)+\tau_K(P)\bigr)
=
\widehat\lambda(P,Q).
\end{aligned}
\]

\smallskip
\noindent
\emph{Case 2: $b=0$.}
In this case
\[
\Delta=ad\in\{\pm1\},
\]
so $a,d\in\{\pm1\}$. Set
\[
\mu_-=\frac{c}{d},
\qquad
\mu_+=0,
\]
and put
\[
x:=p,
\qquad
y:=\mu_-p+q=\frac{c}{d}p+q.
\]
Since $b=0$, the target parameters are
\[
ap+bq=ax
\]
and
\[
\mu_+(ap+bq)+cp+dq
=
cp+dq
=
dy.
\]
Therefore
\[
\widehat\sigma_-(p,q)=U(x)\oplus U(y),
\qquad
\widehat\sigma_+^A(p,q)=U(ax)\oplus U(dy).
\]

Define
\[
\widehat Q_e^0(P,Q):=(S_aP,S_dQ).
\]
Then \eqref{eq:scaling-identity} gives
\[
\begin{aligned}
\widehat Q_e^0\widehat\sigma_-(p,q)(P,Q)
&=
\left(S_aU(x)P,S_dU(y)Q\right)\\
&=
\left(U(ax)S_aP,U(dy)S_dQ\right)\\
&=
\widehat\sigma_+^A(p,q)\widehat Q_e^0(P,Q).
\end{aligned}
\]
The map $\widehat Q_e^0$ is invertible, and again
\[
\widehat Q_e^0\widehat Z=\widehat Z,
\qquad
\widehat\lambda\widehat Q_e^0=\widehat\lambda,
\]
because $S_a$ and $S_d$ fix $K$ and preserve $\tau_K$.

Thus, in either case, we have constructed an invertible map
$\widehat Q_e^0$ satisfying
\begin{equation}
\label{eq:standardized-intertwining}
\widehat Q_e^0\widehat\sigma_-(p,q)
=
\widehat\sigma_+^A(p,q)\widehat Q_e^0,
\qquad
\widehat Q_e^0\widehat Z=\widehat Z,
\qquad
\widehat\lambda\widehat Q_e^0=\widehat\lambda.
\end{equation}

We now return to the original vertex coordinates. Define
\[
Q_e^0:=J_+^{-1}\widehat Q_e^0J_-.
\]
By definition of the standardized actions,
\[
\widehat\sigma_-(p,q)
=
J_-\sigma_{v_-}(\iota_e(h_{p,q}))J_-^{-1}
\]
and
\[
\widehat\sigma_+^A(p,q)
=
J_+\sigma_{v_+}(\iota_{\bar e}(h_{p,q}))J_+^{-1}.
\]
Therefore, using \eqref{eq:standardized-intertwining},
\[
\begin{aligned}
Q_e^0\sigma_{v_-}(\iota_e(h_{p,q}))
&=
J_+^{-1}\widehat Q_e^0J_-
  \sigma_{v_-}(\iota_e(h_{p,q}))\\
&=
J_+^{-1}\widehat Q_e^0
  \widehat\sigma_-(p,q)J_-\\
&=
J_+^{-1}\widehat\sigma_+^A(p,q)
  \widehat Q_e^0J_-\\
&=
\sigma_{v_+}(\iota_{\bar e}(h_{p,q}))
  J_+^{-1}\widehat Q_e^0J_-\\
&=
\sigma_{v_+}(\iota_{\bar e}(h_{p,q}))Q_e^0.
\end{aligned}
\]
This proves the edge-intertwining identity.

Finally, recall that
\[
Z_{a_\pm}^0=J_\pm^{-1}\widehat Z,
\qquad
\lambda_{a_\pm}^0=\widehat\lambda J_\pm.
\]
The pointed identities in
\eqref{eq:standardized-intertwining} therefore imply
\[
\begin{aligned}
Q_e^0Z_{a_-}^0
&=
J_+^{-1}\widehat Q_e^0J_-
  (J_-^{-1}\widehat Z)\\
&=
J_+^{-1}\widehat Q_e^0\widehat Z
=
J_+^{-1}\widehat Z
=
Z_{a_+}^0
\end{aligned}
\]
and
\[
\begin{aligned}
\lambda_{a_+}^0Q_e^0
&=
\widehat\lambda J_+J_+^{-1}
  \widehat Q_e^0J_-\\
&=
\widehat\lambda\widehat Q_e^0J_-
=
\widehat\lambda J_-
=
\lambda_{a_-}^0.
\end{aligned}
\]
This proves all the required identities.

For the reverse orientation, define
\[
Q_{\bar e}^0:=(Q_e^0)^{-1}.
\]
The reverse intertwining and pointed identities follow immediately by
inverting the preceding relations.
\end{proof}

\begin{remark}[Global consistency of the slopes]
The choices \eqref{eq:mu-b-nonzero} and \eqref{eq:mu-b-zero} are made once for
each geometric edge, in its preferred orientation.  Every internal boundary
component corresponds to a unique incidence of a single geometric edge, so each
$\mu_{v,i}$ is assigned exactly once.  
\end{remark}

\subsection{Peripheral detection} We proceed to verify Condition (iii) in \Cref{prop:combination}. For every boundary incidence $(v,i)$, we have defined
\begin{equation}
\label{eq:all-pointed-data}
Z_{v,i}^0:=(K_{v,i},K),
\qquad
\lambda_{v,i}^0(P,Q)
 :=\frac{\tr(PK_{v,i}^{-1})+\tr(QK^{-1})}{12}.
\end{equation}
Then $\lambda_{v,i}^0(Z_{v,i}^0)=1$.

We verify explicitly the invariance required in condition~\textup{(ii)} of
\Cref{prop:combination}.  If $h=c_i^pt_v^q\in H_{v,i}$, then on the two matrix
summands the action of $\sigma_v(h)$ is
\[
\Ad(K_{v,i}^p)\oplus\Ad(K^{\mu_{v,i}p+q}).
\]
It fixes $(K_{v,i},K)$.  Moreover, commutativity of trace and commutativity of
powers of the same matrix give
\[
\tr\bigl(K_{v,i}^pPK_{v,i}^{-p}K_{v,i}^{-1}\bigr)
 =\tr(PK_{v,i}^{-1})
\]
and
\[
\tr\bigl(K^{\mu_{v,i}p+q}QK^{-\mu_{v,i}p-q}K^{-1}\bigr)
 =\tr(QK^{-1}).
\]
Thus
\begin{equation}
\label{eq:lambda-peripheral-invariant}
\lambda_{v,i}^0\sigma_v(h)=\lambda_{v,i}^0
\qquad(h\in H_{v,i}),
\end{equation}
so $H_{v,i}$ preserves $\ker\lambda_{v,i}^0$.  

\begin{lemma}
\label{lem:peripheral-detector}
For every $g\in G_v$,
\begin{equation}
\label{eq:peripheral-detector}
\lambda_{v,i}^0\bigl(\sigma_v(g)Z_{v,i}^0\bigr)=1
\quad\Longleftrightarrow\quad
g\in H_{v,i}.
\end{equation}
\end{lemma}

\begin{proof}
Write $g=wt_v^m$.  Note that $g\in H_{v,i}$ is the same as $w\in \langle c_i\rangle$, where $c_i$ is the group element represented by the corresponding boundary loop.  Conjugate the representation so that the fixed point of \(c_i\) is infinity and write
\[
\rho_v(w)=\begin{pmatrix}r&s\\\gamma&u\end{pmatrix},
\qquad
U_i=\begin{pmatrix}1&\xi\\0&1\end{pmatrix},
\qquad \xi\neq0.
\]
Since $\det\rho_v(w)=1$, direct multiplication gives
\begin{equation}
\label{eq:Fuchsian-trace}
\tr\bigl(\rho_v(w)U_i\rho_v(w)^{-1}U_i^{-1}\bigr)
 =2+\xi^2\gamma^2.
\end{equation}
The factor $\eta_v(w)$ cancels in the tensor commutator.  Recall that $Z_{v,i}^0 = (K_{v,i},K)$, so $
\sigma_v(g)(K_{v,i},K)
=
\bigl(\alpha_v(g)K_{v,i}\alpha_v(g)^{-1},
      \beta_v(g)K\beta_v(g)^{-1}\bigr)$. Consequently,

      \begin{equation}
\label{eq:act-trace}
\begin{aligned}
\lambda^0_{v,i}\bigl(\sigma_v(g)Z^0_{v,i}\bigr)
=\frac1{12}\Bigl(
&\operatorname{tr}\bigl(\alpha_v(g)K_{v,i}\alpha_v(g)^{-1}K_{v,i}^{-1}\bigr)\\
+&\operatorname{tr}\bigl(\beta_v(g)K\beta_v(g)^{-1}K^{-1}\bigr)
\Bigr).
\end{aligned}
\end{equation}
Recall that \[
\alpha_v(w)
=
\bigl(\rho_v(w)\otimes\eta_v(w)\bigr)\oplus I_2,
\qquad
K_{v,i}=(U_i\otimes I_2)\oplus I_2.
\]
Hence

\[
\begin{aligned}
\alpha_v(w)K_{v,i}\alpha_v(w)^{-1}K_{v,i}^{-1}
&=
\Bigl(
(\rho_v(w)\otimes\eta_v(w))
(U_i\otimes I_2)\\
&\hspace{1.5cm}\cdot
(\rho_v(w)^{-1}\otimes\eta_v(w)^{-1})
(U_i^{-1}\otimes I_2)
\Bigr)\oplus I_2\\
&=
\bigl(\rho_v(w)U_i\rho_v(w)^{-1}U_i^{-1}\bigr)
\otimes I_2
\ \oplus I_2.
\end{aligned}
\]
Hence the trace in
the first summand of (\ref{eq:act-trace}) 
\[
\operatorname{tr}\bigl(\alpha_v(g)K_{v,i}\alpha_v(g)^{-1}K_{v,i}^{-1}\bigr) = 2(2+\xi^2\gamma^2)+2=6+2\xi^2\gamma^2.
\]
The trace in the second $\Mat_6$-summand is $6$, because the image of $G_v$
under $\beta_v$ centralizes $K$; see the end of \S \ref{Section:rep-prod-group}.  Therefore
\begin{equation}
\label{eq:detector-value}
\lambda_{v,i}^0\bigl(\sigma_v(g)Z_{v,i}^0\bigr)
 =1+\frac{\xi^2\gamma^2}{6}.
\end{equation}
Over $\R$, this equals $1$ exactly when $\gamma=0$, which says that the
projective element $\bar\rho_v(w)$ fixes the cusp at infinity. But this is the same as saying $g\in\langle c_i,t_v\rangle$. In fact, let
$\Gamma:=\bar\rho_v(\pi_1(Q_v))\leq\PSL_2(\R)$ and let
$p(z)=z+\xi$ be the primitive parabolic corresponding to $c_i$.  Any element
of $\Gamma$ fixing infinity is represented by an affine map $h(z)=az+b$ with
$a>0$.  If $a\neq1$, then
\[
h^nph^{-n}(z)=z+a^n\xi.
\]
For one of the two directions $n\to\pm\infty$, these nontrivial translations
converge to the identity, contradicting discreteness of the Fuchsian group.
Hence $a=1$, so every element of the stabilizer is a translation.  The
translation parameters form a discrete subgroup of $(\R,+)$ and are therefore
cyclic.  Since we assume that $p$ is primitive (see the beginning of \S \ref{Section:rep-prod-group}), it generates that cyclic stabilizer.
Consequently
$\bar\rho_v(w)\in\langle\bar\rho_v(c_i)\rangle$.  Faithfulness of
$\bar\rho_v$ gives $w\in\langle c_i\rangle$.  There is no restriction on $m$,
so this is equivalent to $g\in\langle c_i,t_v\rangle=H_{v,i}$.
\end{proof}

\subsection{Combination and embedding of the coefficient field into
\texorpdfstring{$\C$}{C}}

Let $k_0\subset\R$ be the countable field generated by the finitely many
matrix entries occurring in the vertex representations, pointed data, and
edge intertwiners.  By
\Cref{lem:Theta-faithful,lem:arbitrary-gluing,lem:peripheral-detector} and
\eqref{eq:lambda-peripheral-invariant}, all hypotheses of
\Cref{prop:combination} hold with $n=72$.  After enumerating the finitely many
transcendental variables as $x_1,\dots,x_N$, \Cref{prop:combination} yields an embedding
\begin{equation}
\label{eq:73-field-embedding}
\pi_1(M)\hookrightarrow
\GL_{73}\bigl(k_0(x_1,\dots,x_N)\bigr).
\end{equation}

Since \(k_0\) is countable, we may choose real numbers
\(\zeta_1,\ldots,\zeta_N\in\mathbb R\) that are algebraically
independent over \(k_0\). Indeed, this can be done inductively, since
the relative algebraic closure in \(\mathbb R\) of a countable
subfield of \(\mathbb R\) is countable. The assignment
\(x_i\mapsto\zeta_i\) therefore extends uniquely to an injective field
homomorphism
\[
  k_0(x_1,\ldots,x_N)\longrightarrow\mathbb R.
\]
Applying this field embedding entrywise to the preceding
representation preserves faithfulness. Consequently,
\[
  \pi_1(M)\longrightarrow\operatorname{GL}_{73}(\mathbb R)
\]
is injective.

We have therefore proved the following.

\begin{proposition}[Product-block case]
\label{prop:product-block-case}
Let $M$ be a compact graph manifold all of whose JSJ blocks are products
$Q_v\times S^1$, where each $Q_v$ is compact, connected, orientable, has
nonempty boundary, and has positive genus.  Then
\[
\pi_1(M)\hookrightarrow\GL_{73}(\R).
\]
The gluing maps may be arbitrary elements of $\GL_2(\Z)$ in the peripheral
bases fixed in \eqref{eq:edge-inclusions}.
\end{proposition}

\begin{proof}
    We apply \Cref{prop:combination} to the graph of groups corresponding to the graph manifold. It is enough to note that all the hypotheses hold. For each vertex group, the representation $\sigma_v$ is defined in \eqref{eq:Theta}, and is proved to be injective by Lemma \ref{lem:Theta-faithful}.  For every incidence $a$ at $v$, the vector
 $Z_a^0\in V_v$, and the covector
 $\lambda_a^0\in V_v^*$ are defined in \eqref{eq:unstandardized-pointed-data}. The detector property is proved by Lemma \ref{lem:peripheral-detector}. The condition (iv) follows from Lemma \ref{lem:arbitrary-gluing}.
\end{proof}

Theorem \ref{thm:main} is proved by Proposition \ref{prop:product-block-case}, considering the following fact.

\begin{lemma}
\label{lem:relative-product-cover}
Let $M$ be a closed, connected, orientable, irreducible nongeometric graph
manifold.   Then $M$ has a
finite connected cover
\[
\widehat M\longrightarrow M
\]
such that every Seifert block is a product $Q_v\times S^1$, where
$Q_v$ is compact, connected, orientable, has nonempty boundary, and has genus
at least two. 
\end{lemma}

\begin{proof}
By Wang and Yu \cite[Lemma 0.4]{WY}, the graph manifold $M$ is finitely covered by a graph manifold $\widehat M$ with each piece homeomorphic to $F_i \times S^1$, where $F_i$ is a compact orientable surface and the Euler characteristic $\chi (F_i)<0$. Kapovich--Leeb \cite[Lemma~2.1, especially Steps~1--3]{KapovichLeeb} prove that the graph manifold $\widehat M$ (as a manifold obtained by gluing finitely many Seifert manifolds with hyperbolic base orbifolds) has a further finite cover that has each piece a trivial circle bundle over orientable surfaces of genus $g\geq 2.$ 
\end{proof}

\bibliographystyle{alpha}
\bibliography{references.bib}

@misc{DoubaTs,
    author = {Douba, Sami and Tsouvalas, Konstantinos},
    title = {Matrix entries, unipotents, and linearity of amalgams},
    year = {2026},
    note = {\href{https://arxiv.org/abs/2603.23969}{arXiv:2603.23969}},
    eprint = {2603.23969},
    archivePrefix = {arXiv},
    primaryClass = {math.GR}
}

@book {AFW,
    AUTHOR = {Aschenbrenner, Matthias and Friedl, Stefan and Wilton, Henry},
     TITLE = {3-manifold groups},
    SERIES = {EMS Series of Lectures in Mathematics},
 PUBLISHER = {European Mathematical Society (EMS), Z\"{u}rich},
      YEAR = {2015},
     PAGES = {xiv+215},
      ISBN = {978-3-03719-154-5},
   MRCLASS = {57N10 (20F65 57-02 57M27)},
  MRNUMBER = {3444187},
MRREVIEWER = {Thomas\ Koberda},
       DOI = {10.4171/154},
       URL = {https://doi.org/10.4171/154},
}

@book {Beardon,
    AUTHOR = {Beardon, Alan F.},
     TITLE = {The geometry of discrete groups},
    SERIES = {Graduate Texts in Mathematics},
    VOLUME = {91},
      NOTE = {Corrected reprint of the 1983 original},
 PUBLISHER = {Springer-Verlag, New York},
      YEAR = {1995},
     PAGES = {xii+337},
      ISBN = {0-387-90788-2},
   MRCLASS = {22E40 (11F06 20H15 30F35 57N10)},
  MRNUMBER = {1393195},
}

@book {K3,
     TITLE = {K3---a new problem list in low-dimensional topology},
    SERIES = {Mathematical Surveys and Monographs},
    VOLUME = {295},
    EDITOR = {\.{I}nan\c{c} Baykur, R. and Kirby, Robion C. and Ruberman,
              Daniel},
 PUBLISHER = {American Mathematical Society, Providence, RI},
      YEAR = {[2026] \copyright 2026},
     PAGES = {xi+430},
      ISBN = {978-1-4704-8433-0},
   MRCLASS = {57-02 (57Kxx 57Mxx)},
  MRNUMBER = {5077826},
       DOI = {10.1090/surv/295},
       URL = {https://doi.org/10.1090/surv/295},
}

@article {DavisJanuszkiewicz,
    AUTHOR = {Davis, Michael W. and Januszkiewicz, Tadeusz},
     TITLE = {Right-angled {A}rtin groups are commensurable with
              right-angled {C}oxeter groups},
   JOURNAL = {J. Pure Appl. Algebra},
  FJOURNAL = {Journal of Pure and Applied Algebra},
    VOLUME = {153},
      YEAR = {2000},
    NUMBER = {3},
     PAGES = {229--235},
      ISSN = {0022-4049,1873-1376},
   MRCLASS = {20F36 (20F55 57M07)},
  MRNUMBER = {1783167},
MRREVIEWER = {Thomas\ Brady},
       DOI = {10.1016/S0022-4049(99)00175-9},
       URL = {https://doi.org/10.1016/S0022-4049(99)00175-9},
}

@article {Button,
    AUTHOR = {Button, J. O.},
     TITLE = {A 3-manifold group which is not four dimensional linear},
   JOURNAL = {J. Pure Appl. Algebra},
  FJOURNAL = {Journal of Pure and Applied Algebra},
    VOLUME = {218},
      YEAR = {2014},
    NUMBER = {9},
     PAGES = {1604--1619},
      ISSN = {0022-4049,1873-1376},
   MRCLASS = {57M05 (20H20)},
  MRNUMBER = {3188859},
MRREVIEWER = {Wolfgang\ H.\ Heil},
       DOI = {10.1016/j.jpaa.2014.01.001},
       URL = {https://doi.org/10.1016/j.jpaa.2014.01.001},
}

@article {Boileau,
    AUTHOR = {Boileau, Michel},
     TITLE = {Around 3-manifold groups},
   JOURNAL = {Winter Braids Lect. Notes},
  FJOURNAL = {Winter Braids Lecture Notes},
    VOLUME = {5},
      YEAR = {2018},
    NUMBER = {Winter Braids VIII (Marseille, 2018)},
     PAGES = {Exp. No. 2, 26},
      ISSN = {2426-0312},
   MRCLASS = {57M05 (20F34 57K30)},
  MRNUMBER = {4157114},
}

@inproceedings {AgolICM,
    AUTHOR = {Agol, Ian},
     TITLE = {Virtual properties of 3-manifolds},
 BOOKTITLE = {Proceedings of the {I}nternational {C}ongress of
              {M}athematicians---{S}eoul 2014. {V}ol. 1},
     PAGES = {141--170},
 PUBLISHER = {Kyung Moon Sa, Seoul},
      YEAR = {2014},
      ISBN = {978-89-6105-804-9; 978-89-6105-803-2},
   MRCLASS = {57M50 (57M35 57N10)},
  MRNUMBER = {3728467},
MRREVIEWER = {Sylvain\ Maillot},
}

@article {Liu,
    AUTHOR = {Liu, Yi},
     TITLE = {Virtual cubulation of nonpositively curved graph manifolds},
   JOURNAL = {J. Topol.},
  FJOURNAL = {Journal of Topology},
    VOLUME = {6},
      YEAR = {2013},
    NUMBER = {4},
     PAGES = {793--822},
      ISSN = {1753-8416,1753-8424},
   MRCLASS = {57M05 (20F65 57N10)},
  MRNUMBER = {3145140},
MRREVIEWER = {Stephan\ Rosebrock},
       DOI = {10.1112/jtopol/jtt010},
       URL = {https://doi.org/10.1112/jtopol/jtt010},
}

@article {PWMixed,
    AUTHOR = {Przytycki, Piotr and Wise, Daniel T.},
     TITLE = {Mixed 3-manifolds are virtually special},
   JOURNAL = {J. Amer. Math. Soc.},
  FJOURNAL = {Journal of the American Mathematical Society},
    VOLUME = {31},
      YEAR = {2018},
    NUMBER = {2},
     PAGES = {319--347},
      ISSN = {0894-0347,1088-6834},
   MRCLASS = {57M50 (20F65 57M27)},
  MRNUMBER = {3758147},
MRREVIEWER = {Bruno\ P.\ Zimmermann},
       DOI = {10.1090/jams/886},
       URL = {https://doi.org/10.1090/jams/886},
}

@article {PWBoundary,
    AUTHOR = {Przytycki, Piotr and Wise, Daniel T.},
     TITLE = {Graph manifolds with boundary are virtually special},
   JOURNAL = {J. Topol.},
  FJOURNAL = {Journal of Topology},
    VOLUME = {7},
      YEAR = {2014},
    NUMBER = {2},
     PAGES = {419--435},
      ISSN = {1753-8416,1753-8424},
   MRCLASS = {20F65 (57M27)},
  MRNUMBER = {3217626},
MRREVIEWER = {Ian\ Agol},
       DOI = {10.1112/jtopol/jtt009},
       URL = {https://doi.org/10.1112/jtopol/jtt009},
}

@incollection {Kirby,
    AUTHOR = {Kirby, Robion C.},
     TITLE = {Problems in low-dimensional topology},
 BOOKTITLE = {Geometric topology ({A}thens, {GA}, 1993)},
    SERIES = {AMS/IP Stud. Adv. Math.},
    VOLUME = {2},
     PAGES = {35--473},
    EDITOR = {Kirby, Rob},
 PUBLISHER = {Amer. Math. Soc., Providence, RI},
      YEAR = {1997},
      ISBN = {0-8218-0653-X},
   MRCLASS = {57-02},
  MRNUMBER = {1470751},
       DOI = {10.1090/amsip/002.2/02},
       URL = {https://doi.org/10.1090/amsip/002.2/02},
}

@article{KapovichLeeb,
  author  = {Kapovich, Michael and Leeb, Bernhard},
  title   = {{3}-manifold groups and nonpositive curvature},
  journal = {Geom. Funct. Anal.},
  volume  = {8},
  year    = {1998},
  number  = {5},
  pages   = {841--852},
  doi     = {10.1007/s000390050075}
}

@book{Serre,
  author    = {Serre, Jean-Pierre},
  title     = {Trees},
  series    = {Springer Monographs in Mathematics},
  publisher = {Springer-Verlag},
  address   = {Berlin},
  year      = {2003},
  note      = {Translated from the French original by John Stillwell}
}

@article {WY,
    AUTHOR = {Wang, Shicheng and Yu, Fengchun},
     TITLE = {Graph manifolds with non-empty boundary are covered by surface
              bundles},
   JOURNAL = {Math. Proc. Cambridge Philos. Soc.},
  FJOURNAL = {Mathematical Proceedings of the Cambridge Philosophical
              Society},
    VOLUME = {122},
      YEAR = {1997},
    NUMBER = {3},
     PAGES = {447--455},
      ISSN = {0305-0041,1469-8064},
   MRCLASS = {57N10 (05C10)},
  MRNUMBER = {1466648},
MRREVIEWER = {Hugh\ M.\ Hilden},
       DOI = {10.1017/S0305004197001709},
       URL = {https://doi.org/10.1017/S0305004197001709},
}

@article{Wehrfritz73,
  author  = {Wehrfritz, B. A. F.},
  title   = {Generalized Free Products of Linear Groups},
  journal = {Proceedings of the London Mathematical Society},
  volume  = {s3-27},
  number  = {3},
  pages   = {402--424},
  year    = {1973},
  doi     = {10.1112/plms/s3-27.3.402}
}

@article{Shalen79,
  author  = {Shalen, Peter B.},
  title   = {Linear representations of certain amalgamated products},
  journal = {Journal of Pure and Applied Algebra},
  volume  = {15},
  number  = {2},
  pages   = {187--197},
  year    = {1979},
  doi     = {10.1016/0022-4049(79)90033-1}
}

@article{Mal40,
  author = {Mal'cev, A. I.},
  title = {On isomorphic matrix representations of infinite groups},
  journal = {Rec. Math. [Mat. Sbornik] N.S.},
  volume = {8(50)}, number = {3}, year = {1940}, pages = {405--422}
}

@incollection{Sel60,
  author = {Selberg, Atle},
  title = {On discontinuous groups in higher-dimensional symmetric spaces},
  booktitle = {Contributions to Function Theory (Internat. Colloq. Function Theory, Bombay, 1960)},
  publisher = {Tata Institute of Fundamental Research},
  address = {Bombay}, year = {1960}, pages = {147--164}
}

@article{Tit72,
  author = {Tits, Jacques}, title = {Free subgroups in linear groups},
  journal = {J. Algebra}, volume = {20}, year = {1972}, pages = {250--270},
  doi = {10.1016/0021-8693(72)90058-0}
}

@article{Swn67,
  author = {Swan, Richard G.}, title = {Representations of polycyclic groups},
  journal = {Proc. Amer. Math. Soc.}, volume = {18}, number = {3},
  year = {1967}, pages = {573--574}
}

@article{Big01,
  author = {Bigelow, Stephen J.}, title = {Braid groups are linear},
  journal = {J. Amer. Math. Soc.}, volume = {14}, number = {2},
  year = {2001}, pages = {471--486}, doi = {10.1090/S0894-0347-00-00361-1}
}

@article{Kra02,
  author = {Krammer, Daan}, title = {Braid groups are linear},
  journal = {Ann. of Math. (2)}, volume = {155}, number = {1},
  year = {2002}, pages = {131--156}, doi = {10.2307/3062152}
}

@article{HW08,
  author = {Haglund, Fr{\'e}d{\'e}ric and Wise, Daniel T.},
  title = {Special cube complexes}, journal = {Geom. Funct. Anal.},
  volume = {17}, number = {5}, year = {2008}, pages = {1551--1620},
  doi = {10.1007/s00039-007-0629-4}
}

@article{BN22,
  author = {Berlai, Federico and de la Nuez, Javier},
  title = {Linearity of graph products},
  journal = {Proc. Lond. Math. Soc. (3)}, volume = {124}, number = {5},
  year = {2022}, pages = {587--600}, doi = {10.1112/plms.12434}
}

@incollection{Hem87,
  author = {Hempel, John}, title = {Residual finiteness for 3-manifolds},
  editor = {Gersten, S. M. and Stallings, J. R.},
  booktitle = {Combinatorial Group Theory and Topology (Alta, Utah, 1984)},
  series = {Annals of Mathematics Studies}, volume = {111},
  publisher = {Princeton University Press}, address = {Princeton, NJ},
  year = {1987}, pages = {379--396}
}

@book{AF13,
  author = {Aschenbrenner, Matthias and Friedl, Stefan},
  title = {3-Manifold Groups Are Virtually Residually $p$},
  series = {Memoirs of the American Mathematical Society},
  volume = {225}, note = {No. 1058},
  publisher = {American Mathematical Society}, address = {Providence, RI},
  year = {2013}, doi = {10.1090/S0065-9266-2013-00682-X}
}

@article{Wil17,
  author = {Wilkes, Gareth}, title = {Virtual pro-$p$ properties of 3-manifold groups},
  journal = {J. Group Theory}, volume = {20}, number = {5},
  year = {2017}, pages = {999--1023}, doi = {10.1515/jgth-2016-0060}
}

@article{Wise26,
  author  = {Wise, Daniel T. },
  title   = {Linearity of Graph Manifold Groups
},
  year    = {2026},
  eprint  = {2610.13801},
  archivePrefix = {arXiv}
}

\end{document}